\documentclass[preprint]{elsarticle}
\date{September 23, 2024}

\usepackage{csquotes}
\usepackage{a4wide}
\usepackage[dvipsnames]{xcolor}
\usepackage{graphicx,color,booktabs,listings,caption,subcaption,tikz,amsmath,amssymb,placeins,subfiles,tabto, parskip}

\usepackage[titletoc]{appendix}
\usepackage{algpseudocode}
\usepackage{algorithm}
\usepackage{esvect}
\usepackage{gensymb}
\usepackage{inconsolata}
\usepackage{hyperref}
\usepackage{natbib}
\hypersetup{
    colorlinks=true,
    linkcolor=blue,
    filecolor=magenta,      
    urlcolor=cyan,
    pdftitle={Overleaf Example},
    pdfpagemode=FullScreen,
    }
\usepackage{pdfpages}
\usepackage{arydshln}

\newcommand{\D}{\mathrm{d}}

\journal{Computers \& Fluids}

\begin{document}

\begin{frontmatter}



\title{Reduced Data-Driven Turbulence Closure for Capturing Long-Term Statistics}


\author[inst1]{Rik Hoekstra}
\affiliation[inst1]{organization={Centrum Wiskunde \& Informatica, Scientific Computing group},
            addressline={Science Park 123}, 
            city={Amsterdam},
            postcode={1098 XG}, 
            country={the Netherlands}}
\author[inst1,inst2]{Daan Crommelin}
\author[inst1]{Wouter Edeling\corref{cor1}}
\cortext[cor1]{Corresponding author, e-mail address: wouter.edeling@cwi.nl}

\affiliation[inst2]{organization={Korteweg-de Vries Institute for Mathematics, University of Amsterdam},
            addressline={Science Park 105-107}, 
            city={Amsterdam},
            postcode={1098 XG}, 
            country={the Netherlands}}

\begin{abstract}
    We introduce a simple, stochastic, \textit{a-posteriori}, turbulence closure model based on a reduced subgrid scale term. This subgrid scale term is tailor-made to capture the statistics of a small set of spatially-integrated quantities of interest (QoIs), with only one unresolved scalar time series per QoI. In contrast to other data-driven surrogates the dimension of the ``learning problem" is reduced from an evolving field to one scalar time series per QoI. We use an \textit{a-posteriori}, nudging approach to find the distribution of the scalar series over time. This approach has the advantage of taking the interaction between the solver and the surrogate into account. A stochastic surrogate parametrization is obtained by random sampling from the found distribution for the scalar time series. We compare the new method to an \textit{a-priori} trained convolutional neural network on two-dimensional forced turbulence. Evaluating the new method is computationally much cheaper and gives similar long-term statistics. 
\end{abstract}



\end{frontmatter}

\section{Introduction}
Numerical simulations of turbulent fluid flows pose computational challenges due to the vast range of spatial and temporal scales involved. Resolving all scales in real-world applications, such as climate science, remains unfeasible given current computational capabilities \cite{franzke_stochastic_2016}. Consequently, most simulations only capture a part of the relevant scales (called the resolved scales), and rely on an additional subgrid scale (SGS) term to represent the influence of unresolved scales. In this paper, we examine various models for this SGS term, in the context of large eddy simulations (LES). Approaches to model it are known as closure models or subgrid parameterizations, a term more common in climate science.

Sixty years ago, Smagorinsky introduced the still widely popular empirical Smagorinsky model, which relates the SGS term to the local strain rate of a flow \cite{smagorinsky_general_1963}. Since then, a large variety of subgrid models have appeared based on different approaches, sometimes specialized for specific problems.
Deterministic models such as the Smagorinsky model and its dynamic variant \cite{germano1991dynamic} express the SGS term as a function of the resolved scales. However, given the resolved scales, the unresolved scales might still be in any configuration, leading to stochasticity in the SGS term as specified in the Mori-Zwanzig formalism \cite{franzke_stochastic_2016}. This insight has led to several stochastic SGS models \cite{franzke_stochastic_2016,crommelin2021resampling}, for example based on a simple Ornstein-Uhlenbeck process, but also more sophisticated approaches obeying conservation laws \cite{resseguier_data-driven_2020}, or based on maximum entropy \cite{verkley_maximum_2016}. In these models, the stochastic processes in the SGS term are often fitted to data from a high-fidelity reference simulation.

Advancements in deep learning present new opportunities to learn models for the subgrid scales from data, which falls under the header of scientific machine learning. A recent review of this emerging field can be found in \cite{sanderse_scientific_2024}.
Assume the availability of a reference solution to a flow problem obtained from a high-fidelity solver that captures all relevant scales. We refer to simulations that resolve only a portion of the relevant scales as low-fidelity solvers.
In this context, an important distinction highlighted in \cite{sanderse_scientific_2024} is between \textit{a-priori} and \textit{a-posteriori} learning of the subgrid scale (SGS) term for the low-fidelity solvers. 

The \textit{a-priori} approach learns the parameters in the data-driven SGS model directly from the reference solution. From this reference solution the exact SGS term can be computed and the resolved input features can simply be extracted. The neural network is then trained to approximate a mapping from the resolved scales to the SGS term. This approach has shown to be effective in short-term predictions. However, long-term simulations of low-fidelity solvers coupled to an \textit{a-priori} trained neural network are prone to instabilities and biases \cite{rasp_coupled_2020, frezat_posteriori_2022}. To avoid this, numerous works have focused on adding physical constraints to these models, see e.g.,\ \cite{maulik_sub-grid_2018, guan_learning_2023, yuval_use_2021}.

In the \textit{a-posteriori} approach the interaction between the learned SGS model and the low-fidelity solver is included in the learning process. There are several ways to achieve this, which all require us to run the low-fidelity solver coupled to the reference model during training. A prominent approach is loop-unrolling, where the loss function for the model is computed over multiple steps with the coupled low-fidelity solver \cite{frezat_posteriori_2022, shankar_differentiable_2023,list2022learned}. This requires the low-fidelity solver to be differentiable. An alternative which does not require a differentiable solver is nudging. In \cite{sorensen_non-intrusive_2024} the low-fidelity solver is nudged towards the reference solution. This is accomplished by adding a correction term which linearly relaxes the low-fidelity solution to the reference. Conversely, in \cite{rasp_coupled_2020} a learning approach is described where the high-fidelity solution is nudged to the low-fidelity state. The idea here is to inform low-fidelity behavior from the reference solution, when the latter occupies a similar location in state space as the low-fidelity model.

Deep learning can also be used to create stochastic SGS models. In \cite{guillaumin_stochastic-deep_2021} a data-driven stochastic SGS model is learned by training two neural networks, for the mean and standard deviation of the SGS term, conditional on the coarse-grained solution. In \cite{perezhogin_generative_2023} generative models (generative adversarial networks and variational auto encoders) are trained to predict the SGS term, again conditioned on the current resolved solution. This is a promising approach since these models should be able to learn the full conditional distribution of SGS term, including spatial correlation in randomly generated SGS samples. However, these models are data-hungry and training them is computationally expensive.

The deep learning methods mentioned above, aim to learn a mapping to the SGS term directly. This term is of the same dimension as the low-fidelity solution, which despite the LES filtering operation, still has a large number of degrees of freedom. Whilst computational models can generate an enormous amount of data, it is not uncommon to extract a much lower dimensional set of (spatially integrated) variables from these large models. Think for instance of certain outputs of climate models, which are often global, spatially integrated quantities, e.g.\ the global energy or temperature \cite{zanna_data-driven_2020, rasp2018deep}. If we focus on getting the dynamics of these quantities of interest right, the dimensionality of the learning problem decreases significantly. 

This paper introduces a simple subgrid parameterization based on the reduced subgrid scale term introduced in \cite{edeling_reducing_2020}. This subgrid scale term is tailor-made to capture a small set of quantities of interest (QoIs), with only one unresolved scalar time series per QoI, drastically reducing the dimension of the target of the ``learning problem". We improve upon \cite{edeling_reducing_2020} by using a nudging approach to find a more accurate distribution of these scalar time series, which serve as our new training data. Moreover, \cite{edeling_reducing_2020} only considers the training phase of the reduced SGS models. Here we will focus on the more crucial and difficult predictive phase in which the reduced SGS models are extrapolated in time, using a surrogate model of the time series. We will show that a remarkably simple random noise model can capture the long-term (climate) statistics of the turbulent flow problem under consideration, if given enough training data.
Notable related works are presented in \cite{ephrati_2023_data-driven_Rayleigh_Bernard,ephrati_2023_data-driven_Euler_on_sphere}, which correct LES for known biases in long-term statistical moments using a stochastic parameterization.

 We will denote the aforementioned approach as the ``tau-orthogonal method", named after the orthogonality constraints used in deriving the scalar time series $\tau(t)$ of the reduced subgrid scale term. A comparison is made between this method, a physics-based Smagorinsky model, and a classic deep-learning approach in the form of a convolutional neural network. The Environmental Fluid Dynamics Group at Rice University has done extensive research on using convolutional neural networks (CNNs) as subgrid parameterizations in LES simulations for two-dimensional turbulence. Guan et al.\ \cite{guan_stable_2022} successfully trained a CNN to predict the SGS term in the 2D-DHIT (two-dimensional decaying homogenous isotropic turbulence) equations. They later extended this to forced two-dimensional turbulence, decreasing the required amount of training data by adding physical constraints to the CNN \cite{guan_learning_2023} and leveraging transfer learning for unseen climates \cite{subel_explaining_2023}.
Our CNN subgrid parametrization will use a very similar architecture and \textit{a-priori} learning setup as used in these references.

The forced barotropic vorticity equation is used as test bed, which is widely used to evaluate (data-driven) SGS models  \cite{thuburn_cascades_2014,verkley_maximum_2016,guan_stable_2022, guan_learning_2023, subel_explaining_2023} and can be seen as one of the simplest models in large-scale atmospheric and oceanic modelling when temperature variations are not significant. While this application thus has some relevance, we reiterate that the main purpose of this article is to evaluate the predictive capability of tau-orthogonal method on (spatially-integrated) QoIs, using an as small-and-simple data-driven surrogate as possible, especially compared to a CNN with a multitude of tunable parameters.

This paper is structured as follows. In Section \ref{sec:gov} we describe our governing equations and the discretization method used. Sections \ref{sec:data-driven} and \ref{sec:physics-based} detail the data-driven and physics-based closure models used in our study. We then compare their predictive performance in the subsequent results section, and close with our conclusion in Section \ref{sec:con}.

\section{Governing equations \label{sec:gov}}

We will consider the forced two-dimensional Navier-Stokes equations in barotropic vorticity stream-function formulation:
\begin{align} 
    \frac{\partial \omega}{\partial t}+ J({\psi}, w) & =  \nu \nabla^2 \omega + \mu (F-\omega), \label{eq:HF}\\
    \nabla^2\psi&=\omega.\label{eq:HF2}
\end{align}
\noindent 
Here, $\omega$ is the component of the vorticity perpendicular to the x,y-plane, obtained from the velocity field ${\bf V}=(u, v)^T$ via $\omega={\partial v}/{\partial x}-{\partial u}/{\partial y}$. 
Furthermore, $\psi$ denotes the stream function, related to the velocity components as $u=-\partial\psi/\partial y$, and $v=\partial\psi/\partial x$. The non-linear advection term $J$ is given by:
\begin{equation} \label{eq: J}
J(\psi, w):=\frac{\partial \psi}{\partial x} \frac{\partial w}{\partial y}-\frac{\partial \psi}{\partial y} \frac{\partial w}{\partial x}.
\end{equation}
Following \cite{verkley_maximum_2016} and \cite{edeling_reducing_2020}, we solve these equations on a double periodic domain using values for $\nu$ and $\mu$, which lead to a highly turbulent flow. See \ref{app: problem setup} for a further description of our problem setup.

\subsection{Spectral discretization}

The Fourier spectral method \cite{peyret_vorticity-streamfunction_2002} is used to discretize \eqref{eq:HF}, based on the discrete Fourier expansion 
\begin{equation}
    \omega(x_i,y_j) = \frac{1}{N^2}\sum_{\boldsymbol{k}} \hat{\omega}_{\boldsymbol{k}} e^{i(k_1 x_i + k_2 y_j)}, \text{ for } 0 \leq i,j < N,
\end{equation}
where $\hat{\omega}_{\boldsymbol{k}}$ is the Fourier coefficient for the wave number vector $\boldsymbol{k}=(k_1, k_2)^T$ and $N$ is the number of grid points in each direction of the square grid. Time integration is performed using the standard fourth-order Runge-Kutta scheme. The non-linear term $J$ is computed with a pseudospectral technique, using the ``3/2 rule" to remove aliasing errors, see \ref{app:pseudospectral treatment}. The discrete system is solved using double precision arithmetic.

\subsection{Multiscale decomposition}
We aim to learn a subgrid scale term to improve the accuracy of a computationally cheap low-fidelity (LF) model, using data from a more accurate high-fidelity (HF) model. The LF model uses $65 \times 65$ ($N=65$) Fourier modes, such that it resolves only the largest scales of the problem. We use a HF model with $257 \times 257$ ($N=257$) Fourier modes as ground truth, following \cite{verkley_maximum_2016} and \cite{edeling_reducing_2020}. With this setup the HF model resolves a significantly larger part of the turbulence spectrum than the LF model. Note that it is not feasible to run a fully resolved direct numerical simulation of the problem due to the high Reynolds number, $Re \sim O(10^{6}) $. The HF model uses a time step of $\Delta t = 0.01$ days. Our LF simulations use a 10 times larger time step: $\Delta t = 0.1$. Low fidelity fields are obtained with a spectral filter, simply removing all wave numbers from the HF simulation which are not present in the LF model:

\begin{equation} \label{eq: sharp filter}
    \hat{R}_C \,\hat{x}_{\bf k} =
    \begin{cases}
        \hat{x}_{\mathbf{k}} & \lVert\mathbf{k} \rVert_\infty \leq C \\
        0 & \mathrm{otherwise}
    \end{cases}.
\end{equation}
Thus, $C = \lfloor \frac{N_{LF}}{2} \rfloor=32$ is the cutoff frequency, and $\hat{R}_C$ is a square filter, see Figure \ref{fig:filters_a}. In the following, we will write $\bar{x}=\mathcal{F}^{-1}\left(\hat{R}_C\hat{x}_{\bf k}\right)$, such that $\bar{x}$ is a quantity in the physical domain, filtered by $\hat{R}_C$ in the spectral domain. Applying this filter to \eqref{eq:HF} gives

\begin{align} \label{eq: les1}
    \frac{\partial \bar{\omega}}{\partial t}+ \overline{J({\psi}, \omega)} & =  \nu \nabla^2 \bar{\omega} + \mu (F-\bar{\omega}), \\
    \nabla^2 \bar{\psi}&=\bar{\omega}.
\end{align}
Only the term $\overline{J({\psi}, w)}$ depends on unresolved scales (scales which were filtered away by $\hat{R}_C$), since $\overline{J\left(\psi,\omega\right)}\neq J\left(\bar{\psi},\bar{\omega}\right)$ due to non-linearity. To capture this dependency, we introduce the SGS term:
\begin{equation} \label{eq: subgrid tendency}
    \bar{r} = \overline{J(\psi,w)}-J(\bar{\psi},\bar{\omega}).
\end{equation}
With this exact subgrid scale term, the coarse-grained LES-equations \eqref{eq: les1} can be written as:
\begin{equation}\label{eq: les2}
    \frac{\partial \bar{\omega}}{\partial t}+ J({\bar{\psi}}, \bar{\omega}) =  \nu \nabla^2 \bar{\omega} + \mu (F-\bar{\omega})-\bar{r}.
\end{equation}
Since $\bar{r}$ depends on the unresolved scales we cannot compute it exactly during a low-fidelity simulation. The central problem in closure modelling is finding a good approximation $\tilde{r}$ expressed only in terms of resolved quantities, i.e., \
\begin{equation}
    \tilde{r}(\bar{\psi},\bar{\omega}) \approx \bar{r}(\psi,\omega).
    \label{eq:central_problem}
\end{equation}
This approximation might be interpreted in a broad sense. Most (data-driven) closure modelling efforts attempt to learn an ``ideal" model, where the point-wise difference between $\bar{r}$ and $\tilde{r}$ is small at each time step, at least during the ({\it a-priori}) training phase \cite{beck2019deep,maulik_sub-grid_2018}. Instead, we view \eqref{eq:central_problem} in a weaker sense, and assess the quality of the approximation based on statistics of a set of quantities of interest.

\subsection{Quantities of interest}
We will measure the performance of a surrogate $\tilde{r}$ based on its ability to capture the long-term (climate) statistics of a predefined set of $d$ spatially-integrated QoIs: $\{Q_1,\cdots,Q_d\}$. Here, these are defined as;
\begin{equation}
    Q_i(\bar{\omega},\bar{\psi}) = \left(\frac{1}{2\pi}\right)^2\int\limits_0^{2\pi}\int\limits_0^{2\pi} q_i(\bar{\psi},\bar{\omega}) \D x\D y,\quad i=1,\cdots,d.
    \label{eq:qoi_general}
\end{equation}
Note that each $Q_i$ is a function of our filtered fields, such that it can be computed from either the LF solution or the filtered HF solution. While various QoIs can be cast into this form, common quantities of interest are the energy $E$ and enstrophy $Z$
\begin{align} 
    E = -\frac{1}{2} (\bar{\psi}, \bar{\omega}) \approx -\frac{1}{2} \frac{1}{N_{LF}^4} \sum_{\mathbf{k}} \hat{R}_C\hat{\psi}_\mathbf{k} \hat{\omega}^*_\mathbf{k},\nonumber\\
    Z = \frac{1}{2} (\bar\omega, \bar\omega) \approx \frac{1}{2} \frac{1}{N_{LF}^4} \sum_{\mathbf{k}} \hat{R}_C\hat{\omega}_\mathbf{k} \hat{\omega}^*_\mathbf{k},
    \label{eq:EZ}
\end{align}
where $(f,g) = \iint_{[0,2 \pi] \times [0, 2 \pi]} fg \,\D x \D y / (2\pi)^2$, and $(\cdot)^*$ denotes the complex conjugate. After discretization, this integral is approximated by a vector inner product of Fourier coefficients \cite{edeling_reducing_2020}. Energy and enstrophy are of interest for 2D (geophysical) fluid dynamics models, as the reference solution will remove enstrophy close to the cutoff scale, while energy is backscattered to the large scales. Common SGS models can fail to be properly scale-selective, such that e.g.\ backscatter is not captured, or enstrophy is dissipated at an intermediate scale range, see \cite{thuburn_cascades_2014}. It is therefore interesting to examine $E$ and $Z$ over different scales.

\begin{figure}
    \centering
    \begin{subfigure}[b]{0.317\textwidth}
        \centering
        \includegraphics[width=\textwidth]{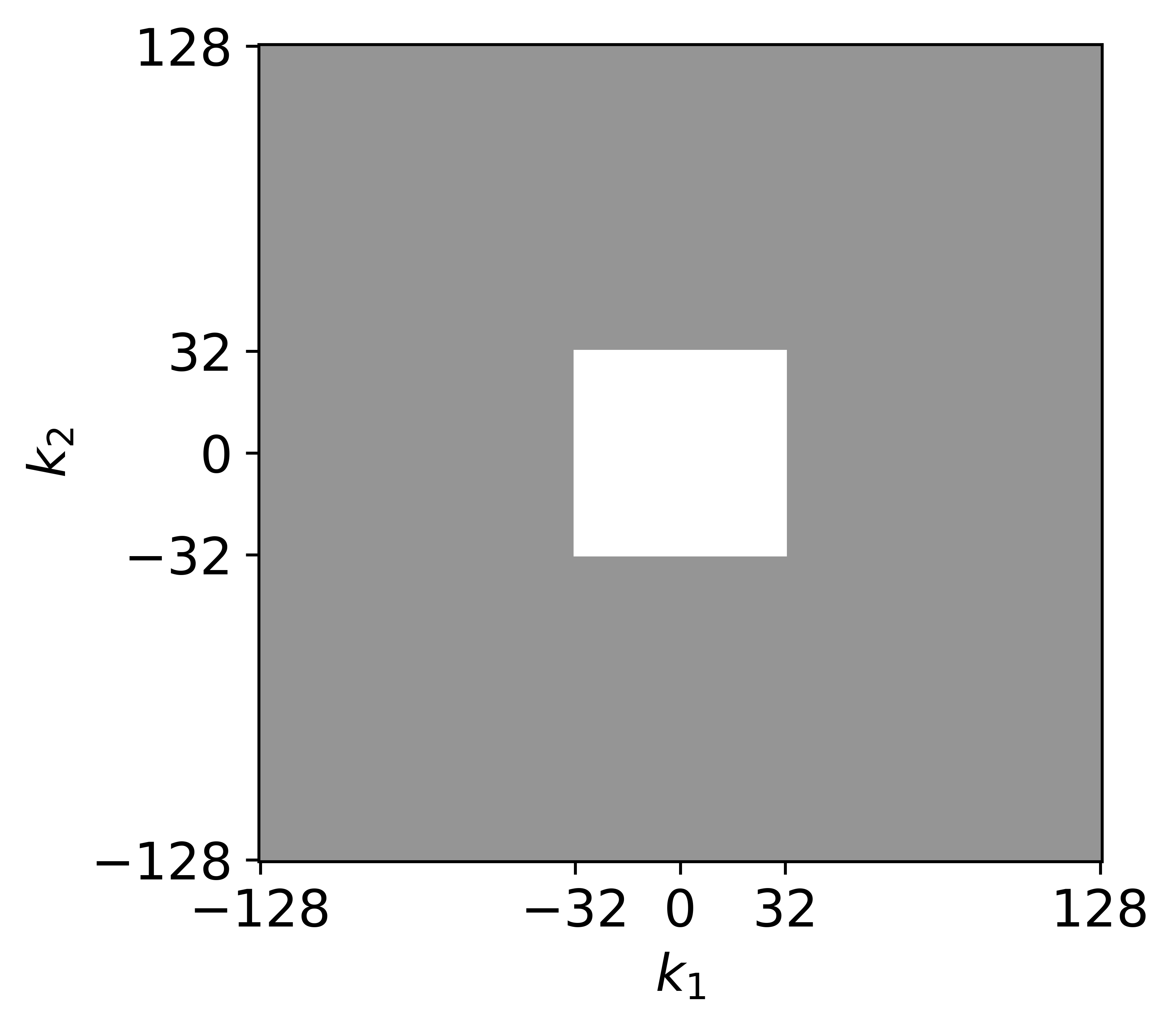}
        \caption{$\hat{R}_C$. \label{fig:filters_a}}
    \end{subfigure}
    \hfill
    \begin{subfigure}[b]{0.29\textwidth}
        \centering
        \includegraphics[width=\textwidth]{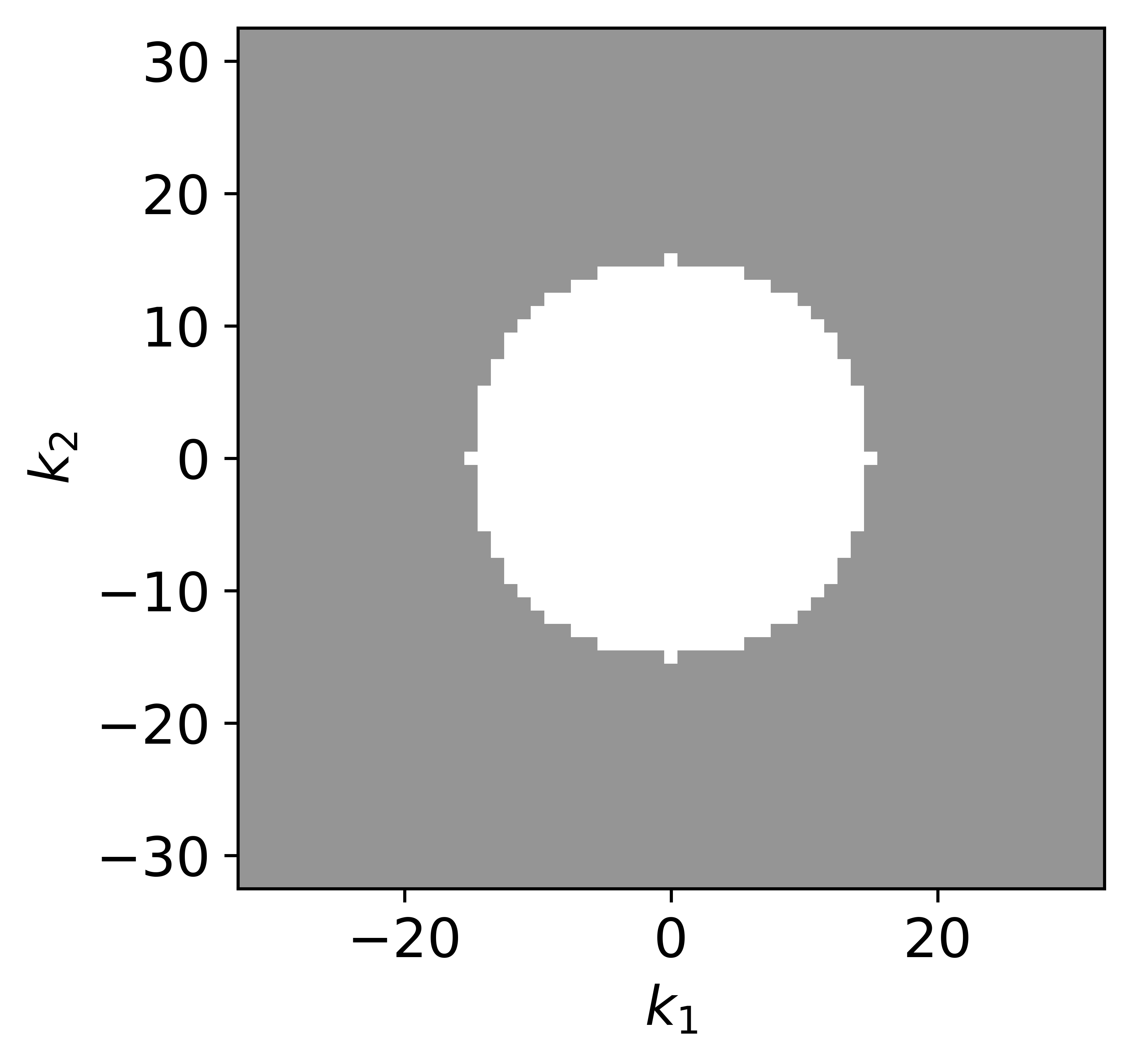}
        \caption{$\hat{R}_{[0,15]}$. \label{fig:filters_b}}
    \end{subfigure}
    \hfill
    \begin{subfigure}[b]{0.29\textwidth}
        \centering
        \includegraphics[width=\textwidth]{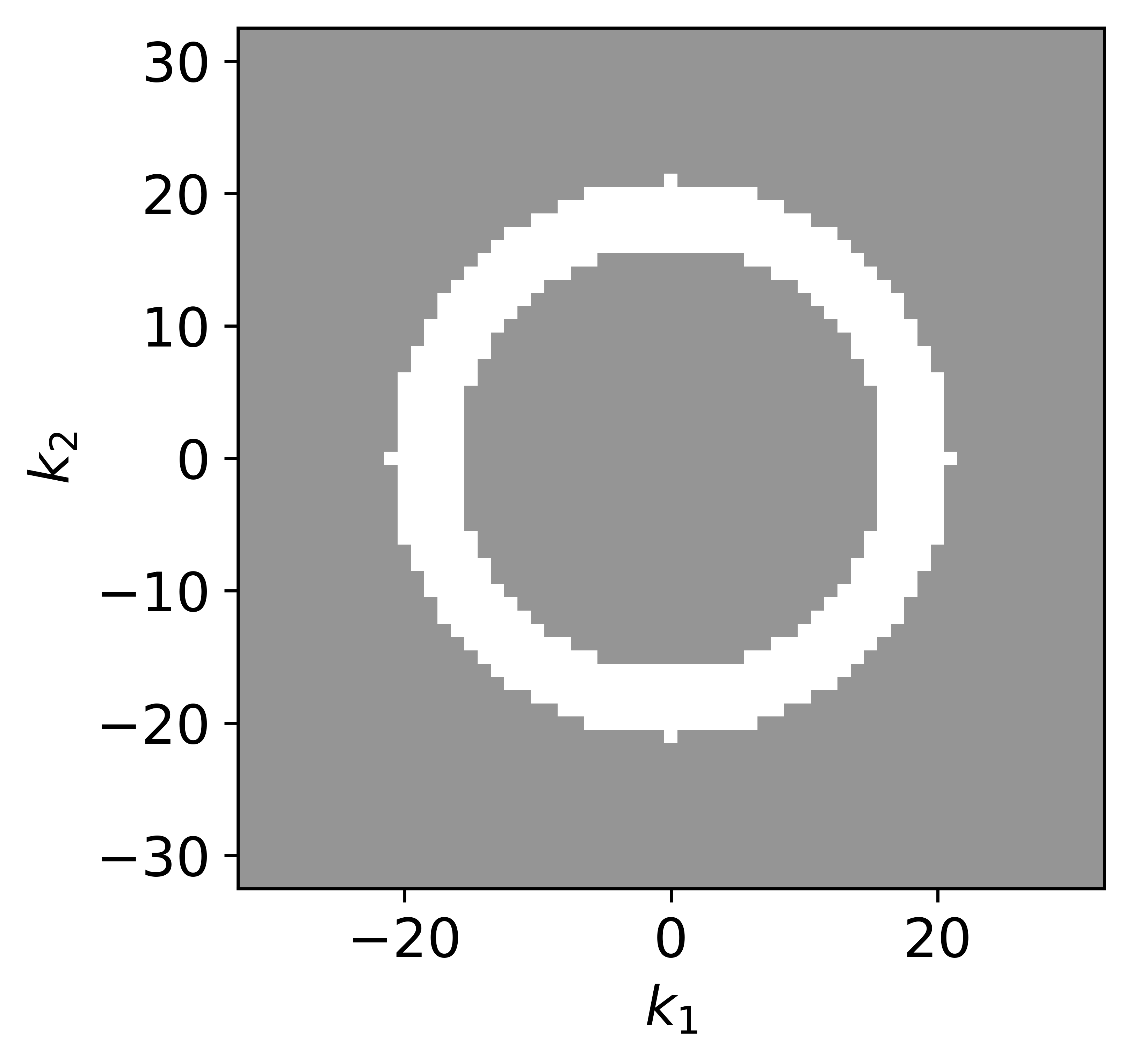}
        \caption{$\hat{R}_{[16,21]}$. \label{fig:filters_c}}
    \end{subfigure}
    \caption{Spectral filters, white = 1, gray = 0.}
    \label{fig:filters}
\end{figure}

To define such ``scale-selective" QoIs, consider the energy contained in a given (one-dimensional) wave number bin $[l,m]$:
\begin{equation}
    E_{[l,m]} = -\frac{1}{2}\frac{1}{N_{LF}^4}\sum_{\mathbf{k}}  \hat{R}_{[l,m]}\hat{\psi}_\mathbf{k} \hat{\omega}^*_\mathbf{k}.
\end{equation}
Here, $[l, m]$ denotes the wave number range $\{{\bf k}\mid l-\frac{1}{2} \leq \sqrt{k_1^2 + k_2^2} < m+\frac{1}{2}\}$ and $\hat{R}_{[l, m]}$ denotes the round filter, defined by
\begin{equation}
    \hat{R}_{[l, m]} \,\hat{x}_{\bf k} =
    \begin{cases}
        \hat{x}_{\mathbf{k}} & ${\bf k}$\in[l, m] \\
        0 & \mathrm{otherwise}
    \end{cases}.
\end{equation}
\noindent
The enstrophy QoIs are defined analogously, and in particular, we will consider the following $d=4$ QoIs:
\begin{equation}
    Q_1 = E_{[0,15]}, Q_2 = Z_{[0,15]}, Q_3 = E_{[16,21]}, Q_4 = Z_{[16,21]},
    \label{eq:4qoi}
\end{equation}
which can be seen as the energy and enstrophy in the large scales and small scales. The corresponding filters are displayed in Figures \ref{fig:filters_b}-\ref{fig:filters_c}. Note that these act on resolved scales only. Hence, for the scale-selective QoI an additional filter is applied on top of $\bar{(\cdot)}$ in the integrand $q_i$, such that for instance
\begin{align}
    E_{[0, 15]} = -\frac{1}{2}\left(R_{[0, 15]}\bar{\psi}, R_{[0, 15]}\bar{\omega}\right),
    \label{eq:E015}
\end{align}
\noindent
where $R_{[0, 15]}$ acts in physical space, such that  $\hat{R}_{[0, 15]}\hat{\psi}_{\bf k}=\hat{R}_{[0, 15]}\hat{R}_C\hat{\psi}_{\bf k}$. For brevity, we will also denote the scale-selective filter for the $i$-th QoI by $R_i$, so\ $Q_1 = -1/2\left(R_1\bar\psi, R_1\bar\omega\right)$.

\section{Data-driven closure models \label{sec:data-driven}}

In this section we describe two different data-driven methods for learning $\tilde{r}$. First the tau-orthogonal method is described. As mentioned, a so-called reduced SGS term is derived, that only focuses on the predefined QoIs.

In Section \ref{section: CNN}, we describe a traditional {\it a-priori} approach where a Convolutional Neural Network (CNN) is trained on full-field  $\overline{r}$ snapshots. This leads to a general surrogate $\tilde{r}$ that closely resembles the exact SGS term $\bar{r}$. It is a black-box approach that must approximate the SGS field that is unresolved at every grid point. 

\subsection{Tau-orthogonal method} \label{section: TO method}

The tau-orthogonal method is based on the following ansatz for the functional form of $\tilde{r}$:
\begin{equation} \label{eq: proposed surrogate}
    \tilde{r}(x, y, t) = \sum_{j=1}^{d} \tau_j(t)P_j(x, y, t),
\end{equation}
\noindent
where $d$ is the number of QoIs, $P_j$ are resolved full-field patterns and $\tau_j$ are unresolved time series, all described in detail below.

The overarching idea behind \eqref{eq: proposed surrogate} is to sacrifice generality for low-dimensionality. By making the concession of not even attempting to learn a perfect surrogate of $\bar{r}$, we gain a reduction in the unresolved degrees of freedom by a factor of $N^2_{LF}/d$ (more than 1000 in our case), as only the $\tau_i$ will be unresolved. This in turn allows us to close the system at the level of the ordinary differential equations (ODEs) that govern our QoIs, instead of at the (high-dimensional) partial differential equation \eqref{eq: les2}. For the scale-selective QoIs these ODEs are given by (see e.g.\ \eqref{eq:E015}):
\begin{equation} \label{eq: ODE q}
    \frac{\mathrm{d}Q_i}{\mathrm{d}t} = \left(\frac{\partial q_i}{\partial R_i\bar\omega}, \frac{\partial R_i\bar\omega}{\partial t}\right) + 
    \left(\frac{\partial q_i}{\partial R_i\bar\psi}, \frac{\partial R_i\bar\psi}{\partial t}\right),
   \quad i=1,\cdots,d,
\end{equation}
which we can rewrite in terms of $\bar{\omega}$ as 
\begin{equation} 
    \frac{\mathrm{d}Q_i}{\mathrm{d}t} = \left(V_i, \frac{\partial R_i\bar\omega}{\partial t}\right),
   \quad i=1,\cdots,d,
   \label{eq:V_i1}
\end{equation}
using;
\begin{equation}
    V_1 = -R_1 \bar{\psi}, \quad V_2 = R_2 \bar{\omega}, \quad V_3 = -R_3 \bar{\psi}, \quad V_4 = R_4 \bar{\omega}.
    \label{eq:V_i2}
\end{equation}
\noindent
For the derivations behind \eqref{eq:V_i1}-\eqref{eq:V_i2}, see e.g.\ \cite{verkley_maximum_2016,edeling_reducing_2020}. Substitution of the LES equations in \eqref{eq:V_i1} gives 
\begin{equation}
    \frac{\mathrm{d} Q_i}{\mathrm{d} t} = \left(V_i, 
    R_i\left[- J({\bar{\psi}}, \bar{\omega}) + \nu \nabla^2 \bar{\omega} + \mu (F-\bar{\omega})-\tilde{r}\right] \right).
    \label{eq: ODE q4}
\end{equation}
\noindent
We lump all resolved-scale contributions into ${\mathrm{d} Q^r}/{\mathrm{d} t}$, and denote the unresolved-scale contribution by ${\mathrm{d} Q^u}/{\mathrm{d} t}$, such that after inserting \eqref{eq: proposed surrogate} into \eqref{eq: ODE q4} we obtain
\begin{equation}
    \frac{\mathrm{d} Q_i}{\mathrm{d} t} = \frac{\mathrm{d} Q_i^r}{\mathrm{d} t} + \frac{\mathrm{d} Q_i^u}{\mathrm{d} t} =  
    \frac{\mathrm{d} Q_i^r}{\mathrm{d} t} - \sum_{j=1}^{d}\tau_j(t) \left(V_i, R_i P_j \right).
    \label{eq: ODE q3}
\end{equation}
\noindent
The last term in \eqref{eq: ODE q3} is the SGS term at the ODE level resulting from the assumption made on $\tilde{r}$ at the PDE level.

The closure problem in \eqref{eq: ODE q3} consist of i) specifying the unresolved $\tau_j$, and ii) specifying the resolved patterns $P_j$. The latter are chosen to simplify the construction of the former. In particular, the full-field patterns $P_j$ are constructed such that they only influence one QoI each. The following orthogonality condition imposes this
\begin{equation} \label{eq: orthogonality constr}
    \left( V_i, R_i P_j \right) = 0 \text{ if } i \neq j.
\end{equation}
\noindent
This simplifies the equations to
\begin{equation}
    \frac{\mathrm{d} Q_i}{\mathrm{d} t} =  \frac{\mathrm{d} Q_i^r}{\mathrm{d} t} - \tau_i(t) \left(V_i, R_i P_i \right),
\end{equation}
\noindent
where the $\tau_i$ are related to the unresolved scale contributions via
\begin{equation}\label{eq: tau in terms of dq}
      \tau_i =- { \frac{\mathrm{d} Q_i^u}{\mathrm{d} t} } \bigg/ {\left(V_i, R_i P_i \right)}.
\end{equation}
\noindent
To fulfil the orthogonality constraint \eqref{eq: orthogonality constr}, the full-field patterns are defined as
\begin{equation}
    P_i(x, y, t) = \sum_{j=1}^d c_{i,j}(t) T_{j}(\bar{\psi}, \bar{\omega}).
    \label{eq:Pi}
\end{equation}
\noindent
The $T_{j}$ are user-specified resolved basis functions, which are a function of the resolved solution $\{\bar{\psi}, \bar{\omega}\}$. For simplicity, we choose the same basis for $P_j$ to which they must be made orthogonal, i.e.\ $T_{j} = V_j$. While other choices could be feasible, we previously obtained good and numerically stable results with this particular choice \cite{edeling_reducing_2020}. For every $P_i$ there are $d-1$ orthogonality constraints and $d$ coefficients $c_{ij}$. Hence, if we set $c_{ii}(t)=1$ and insert \eqref{eq:Pi} into \eqref{eq: orthogonality constr}, we obtain $d$ closed linear systems for the remaining $c_{ij}(t)$ coefficients.

As an example, for $i=1$ and $d=4$ this would look like:
\begin{align}
P_1 &= T_1 + c_{12}T_2 + c_{13}T_3 + c_{14}T_4
 \quad\bot \quad V_2, V_3,V_4,
\end{align}
\noindent
which via \eqref{eq: orthogonality constr} gives a linear system
\begin{align}
    0=\left(V_j, R_j T_1\right) +
    \sum_{k\in\{1,\cdots,4\}\backslash\{1\}} c_{1k} \left(V_j, R_j T_k\right),
    \quad j=2,3,4,
    \label{eq:c_1j}
\end{align}
of size $3\times 3$ that can be solved for the $c_{1j}(t)$ at every time step, and likewise for $i=2,\cdots,4$. 

Note that since our $q_k$ are filtered through $R_k$, the term $R_j T_k=R_jR_kV_k$ appearing in \eqref{eq:c_1j} is zero if $R_j$ and $R_k$ are disjoint, i.e.\ if $\hat{R}_j$ and $\hat{R}_k$ have no overlap in spectral space, which is the case for our filters. We can use this to simplify the analytic expression that underpins $\tilde{r}$, see \ref{app:analytic} for the derivation of the functional form of our surrogate. That said, this analytic expression is not explicitly used in practice, as we numerically solve for the $c_{ij}$ coefficient via linear systems such as \eqref{eq:c_1j}.

\subsubsection{Nudging approach}

We expressed our reduced surrogate term \eqref{eq: proposed surrogate} via the relations in \eqref{eq: tau in terms of dq} and \eqref{eq:Pi} as a function of the resolved solution fields $\{\bar{\psi}, \bar{\omega}\}$ and the unresolved-scale contribution to the ODEs of the QoIs:
\begin{equation}
    \tilde{r}\left(\bar{\psi}(t),\bar{\omega}(t)\right) = \tilde{r}\left(\frac{\mathrm{d} Q^u}{\mathrm{d} t}\bigg|_{t}, \bar{\psi}(t), \bar{\omega}(t)\right).
\end{equation}
\noindent
To nudge the QoI trajectories in the LF solution towards the trajectories in the HF solution, we discretize this with a predictor-corrector like scheme. Here, the predictor step comprises a step of the LF solver without SGS term and the correction is calculated based on the predicted values. Let us denote the predicted fields with an extra asterisk in the super-script: i.e. $\{\bar{\omega}^{n+1^*}\}$. We approximate the unresolved-scale contribution to the ODE as
\begin{equation}
    \frac{\D Q_i^u}{\mathrm{d} t}\bigg|_{t_{n+1}} \approx \frac{\Delta Q_i^u(t_{n+1})}{\Delta t} := \frac{Q_i(\bar{\psi}^{n+1}_{HF},\bar{\omega}^{n+1}_{HF}) - Q_i(\bar{\psi}^{n+1^*},\bar{\omega}^{n+1^*})}{\Delta t},
    \label{eq:dQu}
\end{equation}
where $Q_i(\bar\psi^{n+1}_{HF},\bar\omega^{n+1}_{HF})$ is the value of the reference QoI at time $t_{n+1}$, computed with the HF vorticity and stream function projected on the LF grid. Note that we assume that the previous corrector step was successful, i.e.\ that $\Delta Q_i^u(t_n)=0$. The correction step is given by
\begin{align}
    \bar{\omega}^{n+1} & = \bar{\omega}^{n+1^*} + \int_{t_n}^{t_{n+1}} \tilde{r}\left(\frac{\Delta Q_i^u(t_{n+1})}{\Delta t}, \bar{\psi}^{n+1^*}, \bar{\omega}^{n+1^*}\right) \D t, \\
                       & = \bar{\omega}^{n+1^*} + \tilde{r}\left(\Delta Q_i^u(t_{n+1}), \bar{\psi}^{n+1^*}, \bar{\omega}^{n+1^*}\right),
\end{align}
where we can move the $1/\Delta t$ factor because $\tilde{r}$ is linear in ${\mathrm{d} Q^u}/{\mathrm{d} t}$. In the following we will use $\Delta Q_i$ to denote the difference between the value of the $i$-th QoI in the reference simulation and the current (predicted) low fidelity state, where we dropped the superscript $u$, to ease the notation.

Inside the training domain, i.e. when HF reference trajectories are available for the QoIs, this approach tracks these trajectories with high accuracy, see Section \ref{sec:TO_training}. In \ref{app:error in correction} we derive a bound on the difference between the HF reference trajectories and the LF trajectories at time $t$ with this correction scheme.

\subsection{CNN} \label{section: CNN}

We use the CNN architecture tuned in \cite{guan_stable_2022}, and applied to essentially the same governing equations as studied here. It contains 10 convolutional layers with convolutional depth 64 and filter size $5 \times 5$. ReLU activation functions are used, except for the linear output layer. It is evaluated using single precision arithmetic, as is common practice in machine learning applications. The only difference is that our layers use periodic padding, in line with our periodic boundary conditions, whereas \cite{guan_stable_2022} used zero padding. The CNN is trained to approximate the mapping
\begin{equation}
    \left\{ \bar{\psi}/\sigma_{\bar{\psi}}, \bar{\omega}/\sigma_{\bar{\omega}} \right\} \to \left\{ \bar{r}/\sigma_{\bar{r}} \right\},
    \label{eq:cnn}
\end{equation}
where $\sigma_{\bar\psi}$, $\sigma_{\bar\omega}$, $\sigma_{\bar{r}}$ denote the standard deviations of $\bar\psi$, $\bar\omega$ and $\bar{r}$ respectively, calculated over all training samples. As such, \eqref{eq:cnn} represents a mapping from $\mathbb{R}^{2 \times N_{LF} \times N_{LF}}$ to $\mathbb{R}^{N_{LF} \times N_{LF}}$.

For training, the AdamW optimizer \cite{loshchilov_decoupled_2019,kingma_adam_2017} is used to minimize the mean-squared error between the predicted and true SGS term. The CNNs are trained for 100 epochs, with a learning rate of $2 \cdot 10^{-4}$ during the first 50 epochs and $1 \cdot 10^{-4}$ for the last 50 epochs.

\begin{figure}[htp]
    \captionsetup{justification=centering}
     \centering
     \begin{subfigure}[b]{0.23\textwidth}
         \centering
         \includegraphics[width=\textwidth]{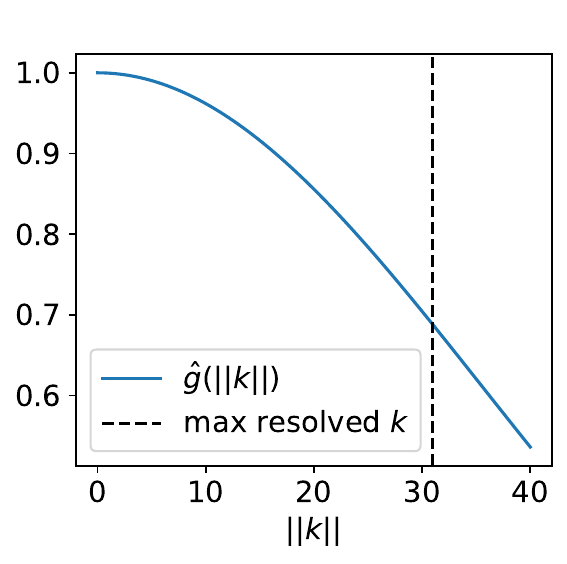}
         \caption{Gaussian filter in Fourier space.}
         \label{fig: Gaussian filter}
     \end{subfigure}
     \hfill
     \begin{subfigure}[b]{0.76\textwidth}
         \centering
         \includegraphics[width=\textwidth]{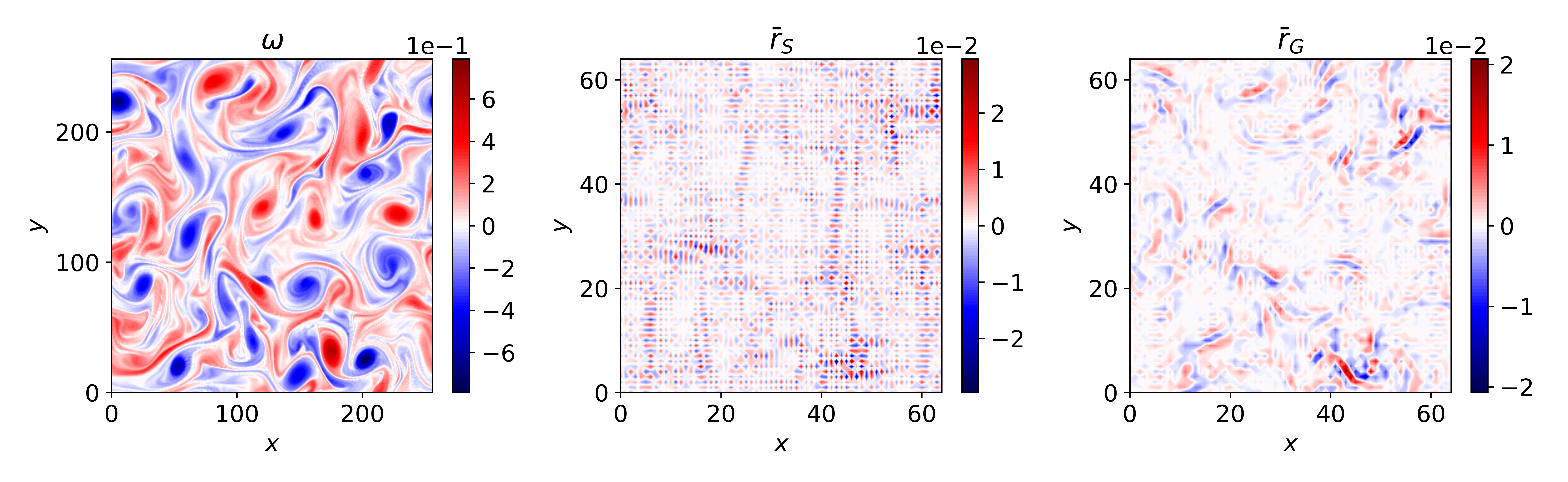}
        \caption{HF snapshot and corresponding SGS snapshots for the sharp spectral and Gaussian filters.
        }
        \label{fig: sgs term for filters}
     \end{subfigure}
     \caption{The Gaussian filter and various flow snapshots.}
\end{figure}

In \cite{guan_stable_2022, guan_learning_2023, subel_explaining_2023} the LF simulations are derived from the HF simulations using a Gaussian filter. To coarse-grain the HF solution to the LF grid in our setup, \eqref{eq: sharp filter} is used. This sharp spectral filter is compact and easy to work with in Fourier space. Yet in physical space, where the CNN learns the subgrid term, the sharp spectral filter is non-local. This represents a challenge to the CNN since its kernels are local. Hence the use of the Gaussian filter, which is compact in both physical and Fourier space \cite{pope_turbulent_2000}. By applying the Gaussian filter after coarse-graining, the overall operation is local. The Gaussian filter in Fourier space, applied element-wise to all Fourier coefficients, is defined as
\begin{equation}
    \hat{G}\hat{u}_{\boldsymbol{k}} = 
        e^{- \lVert \boldsymbol{k} \rVert_2^2 \Delta_{F}^2/24} \hat{u}_{\boldsymbol{k}},
\end{equation}
where $\Delta_{F} = 2 \Delta_{LF} = 2 \cdot \frac{2 \pi}{N_{LF}}$ is the filter width. Figure \ref{fig: Gaussian filter} contains a plot of the Gaussian kernel in Fourier space, for $N_{LF} = 65$.

A different filter in the coarse-graining operation leads to a different SGS term. Figure \ref{fig: sgs term for filters} shows a HF solution field for $\omega$, the corresponding SGS term $\bar{r}_S$ for the sharp spectral filter, and the SGS term when the Gaussian filter is used on top of the sharp filter $\bar{r}_G$. We will see in Section \ref{sec:results} that after training, a CNN is much better capable of predicting $\bar{r}_G$ than $\bar{r}_S$. 

\section{Physics-based closure model \label{sec:physics-based}}

The Smagorinsky model \cite{smagorinsky_general_1963} is a commonly used physics-based (eddy-viscosity) subgrid parameterization, which assumes the subgrid tensor is linearly proportional to the coarse strain-rate tensor $\bar{S}_{ij} = \frac{1}{2}\left(\frac{\partial \bar{u}_i}{\partial x_j}+\frac{\partial \bar{u}_j}{\partial x_i}\right)$. It replaces the viscous term in the filtered Navier-Stokes equations by
\begin{equation} \label{eq: smag velosity}
    \nabla \cdot \big( 2(\nu - \nu_s(\boldsymbol{x}))  \bar{S}_{ij}(\boldsymbol{x})\big),
\end{equation}
\noindent
where $\nu_s$ is given by
\begin{equation}
    \nu_s = (C_s \Delta)^2 \left( 2 \bar{S}_{ij} \bar{S}_{ji} \right) ^{1/2}.
\end{equation}
Here, $C_s$ is a tunable constant, $\Delta$ the subgrid scale characteristic length, and repeated tensor indices are used for Einstein summation. The Smagorinsky constant can be approximated for isotropic turbulence by \cite{lesieur_large-eddy_2005, lilly_representation_1967}:
\begin{equation}
        C_s \approx \frac{1}{\pi} \left( \frac{3 C_k}{2}\right)^{-\frac{3}{4}}.
\end{equation}
Which yields $C_s \approx 0.18$ for a Kolmogorov constant equal to $C_k = 1.4$. However, the optimal value of $C_s$ is problem dependent. For general problems, a value of 0.1 is a common choice \cite{lesieur_large-eddy_2005}. In literature a large range of values is used, from $0.0065$ in the bulk of channel flows \cite{gatski_simulation_1996} to $1.0$ when fitted offline to the subgrid scales in 2D decaying turbulence \cite{maulik_sub-grid_2018}. \cite{shankar_differentiable_2023} found a value of $0.172$ using gradient descent in a two-dimensional turbulence setup, similar to the one studied here.
\noindent
Applying the curl operator to \eqref{eq: smag velosity} gives us the viscous term in the vorticity formulation
\begin{equation}
    \nu \nabla^2 \bar{\omega} - \nabla \times [ 2 \nabla \cdot \left( \nu_s(\boldsymbol{x}) \bar{S}_{ij}(\boldsymbol{x}) \right)].
\end{equation}
Where $\nu_s$ is given by \cite{guan_stable_2022}
\begin{equation}
    \nu_s(\bar{\omega}, \bar{\psi}) = (C_s \Delta)^2  \sqrt{4 \left( \frac{\partial^2 \bar{\psi}}{\partial x \partial y} \right)^2 + \left( \frac{\partial^2 \bar{\psi}}{\partial x ^2} -  \frac{\partial^2 \bar{\psi}}{\partial y ^2}\right)^2}.
\end{equation}
The subgrid scale surrogate is then given by
\begin{equation} \label{eq: smag sgs tendency}
    \tilde{r}_{Smag}(\bar{\omega}, \bar{\psi})  =  2 \nabla \times [ \nabla \cdot \left( \nu_s(\bar{\omega}, \bar{\psi}) \bar{S}_{ij}(\bar{\omega}, \bar{\psi}) \right)].
\end{equation}

\section{Results \label{sec:results}}

In this section we will briefly discuss the training results of the tau-orthogonal method, and examine the amount of training data it requires. Subsequently we will compare the predictions of the tau-orthogonal method, the CNN and the physics-based Smagorinsky model.

\subsection{Tau-orthogonal results}

\subsubsection{Training error \label{sec:TO_training}}

\begin{figure}[h!]
    \centering
    \includegraphics[width = \textwidth]{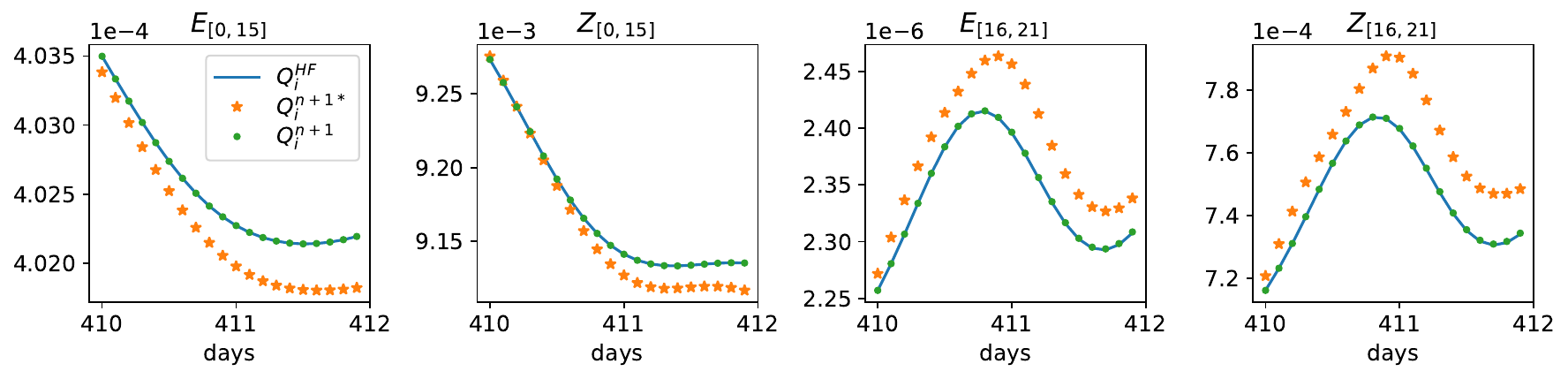}
    \caption{Tracking QoIs with the predictor-corrector approach.}
    \label{fig: TO tracking in training}
\end{figure}

In Figure \ref{fig: TO tracking in training} we demonstrate the tracking procedure in the training phase, where the $\Delta Q_i$ are computed with respect to an HF reference simulation. We plotted the QoI values after the prediction step (denoted by the *) and after the correction step. This shows that, as far as QoI prediction is concerned, there is no discernible difference between using the exact SGS term $\overline{r}$ or the reduced surrogate $\tilde{r}$ during training. For more training results we refer to \cite{edeling_reducing_2020}.

\subsubsection{Convergence of distributions}

Hence, during the training phase we faithfully reproduce the HF $Q_{i}$ trajectories. However, our main focus here is the predictive phase, in which we replace $\Delta Q_i$ with a stochastic data-driven surrogate in order to capture the HF $Q_{i}$ long-term distributions. In this section we therefore compare the convergence of the (marginal) $\Delta Q_i$ and HF $Q_{i}$ distributions.

To compare one-dimensional (empirical) distributions we use the Kolmogorov-Smirnov (KS) test statistic, which measures the distance between two cumulative density functions (CDFs) $F(x)$ and $G(x)$ via
\begin{equation}
    KS\Big(F(x),G(x)\Big) = \max_x \Big\lvert F(x)-G(x) \Big\rvert.
    \label{eq:KS}
\end{equation}

Instead of using a pointwise measure as \eqref{eq:KS}, it might be appropriate to use a similarity measure that integrates over the support of the two pdfs, e.g.\ the Hellinger or Wasserstein distance. However, to be able to compare the distances obtained from different QoIs binning or normalization choices have to be made, because the distributions of the QoI live on different scales. We opted to retain the simpler KS statistic, avoiding these choices.
We found our results to be robust under a change to the Hellinger distance (where we binned based on the spread of the reference distribution).
 To get an intuition of what can be considered a small KS-distance, we compare the marginals of our QoIs over the first 10000 days of a 20000-day HF simulation to the distribution over the full 20000 days, see Table \ref{tab: KS HF converged}. 
\begin{table}[htp]
    \centering
        \begin{tabular}{c|c}
             QoI & KS-distance \\ \toprule
             $E_{[0,15]}$& 0.029\\
             $Z_{[0,15]}$& 0.023\\
             $E_{[16,21]}$& 0.009\\
             $Z_{[16,21]}$& 0.008\\\bottomrule
        \end{tabular} 
        \caption{Similarity in distribution of QoIs in first half of a 20000-day HF simulation and the full simulation.}
        \label{tab: KS HF converged}
\end{table}
\begin{figure}[h!]
    \centering
    \includegraphics[width = 0.6 \textwidth]{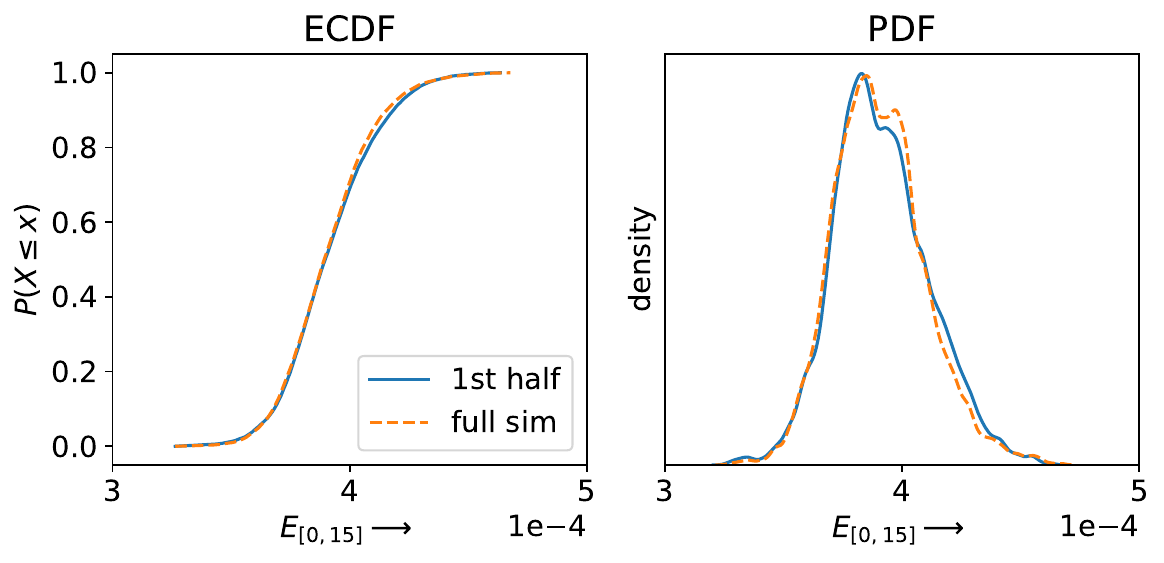}
    \caption{Empirical cumulative distribution functions and kernel-based probability density functions for the marginal distribution of $E_{[0,15]}$ for the first 10000 days and full 20000 days of the HF simulation. }
    \label{fig: ecdf_pdf_z_015 HF}
\end{figure}

\begin{figure}
    \centering
    \includegraphics[width=0.9\textwidth]{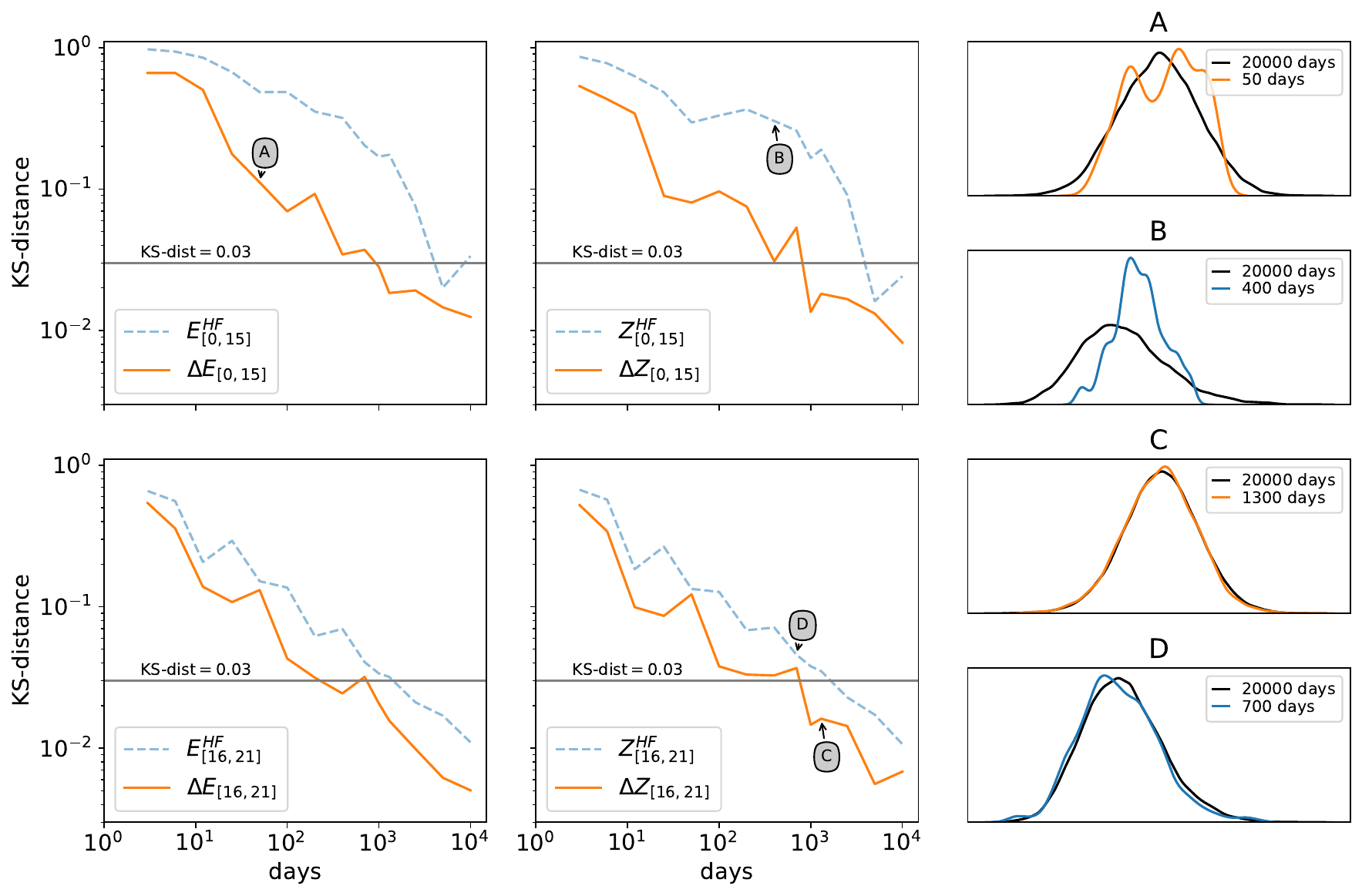}
    \caption{Convergence of marginals for $\Delta Q_i$ and HF $Q_i$. The large panels show the KS-distance between the marginal in the first $x$ days and the reference marginal over the full 20000 days. A KS distance of 0.03 corresponds to the largest value from Table \ref{tab: KS HF converged}, below which we consider the PDF to be converged. The small panels show some PDFs for points in the convergence graph, with the marginal from the full 20000 days in black.}
    \label{fig: convergence of dQ and Q_HF}
\end{figure}

We first inspect the convergence of the distributions of $Q_{i, HF}$ and $\Delta Q_i$ in terms of the KS-distance, see Figure \ref{fig: convergence of dQ and Q_HF}, where we compare the distribution over the first $x$ days of a simulation to the distribution over the full 20000 simulation days.\footnote{We exclude the first 200 days as burn-in period for the low fidelity simulation.} Especially the corrections on the large scales converge faster than the corresponding QoI distributions in the HF simulation, resulting in convergence of all correction distributions within 1000 days, which is considerably faster than the 10000 days needed for the QoI distributions in the HF solver. The probability density functions (PDFs) of the corrections $\Delta Q_i$ (A \& C) resemble Gaussian distributions, which we will try to exploit later.

\begin{figure}
    \centering
    \includegraphics[width=0.9\linewidth]{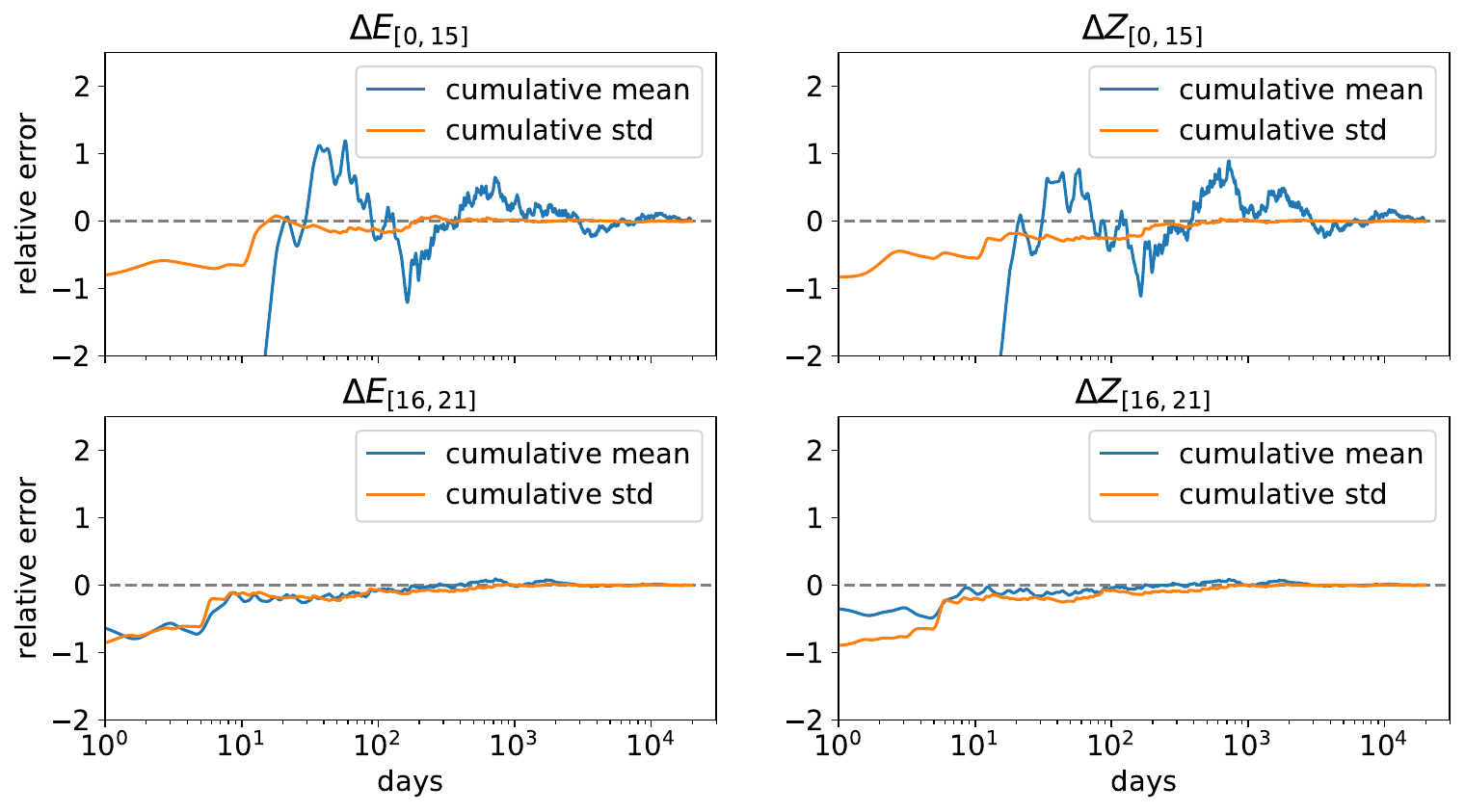}
    \caption{Convergence of the cumulative mean and the cumulative standard deviation of the marginals for $\Delta Q_i$.}
    \label{fig: convergence of mean dQ}
\end{figure}

In Figure \ref{fig: convergence of mean dQ}, we take a closer look at the convergence of the relative error in the mean and the standard deviation of the $\Delta Q_i$. The relative error in the standard deviation quickly decreases. However, for the large-scale QoI ($E_{[0,15]}$ and $Z_{[0,15]}$) the relative error in the mean remains larger, because the long term means are small. Getting these means correct will turn out to be of crucial importance for the performance of the tau-orthogonal method.


\subsubsection{Tau-orthogonal prediction results}

TO surrogate construction was already simplified by the dimension reduction, as only $\Delta Q_i(t)$ must be learned from data. It is of course still possible to fit a (large) machine-learning model to this time-series data. However, in the following we will impose a very simple functional form on the $\Delta Q_i$ surrogates. In particular, we will examine the following 3 simple data-driven $\Delta Q_i$ surrogates. All are solely focused on reproducing the distribution of $\Delta Q_i$. As such, no correlation between time steps is present, we are replacing $\Delta Q_i$ with a data-driven noise model.
\begin{enumerate}
    \item {\bf Plain resampling}: here we simply resample $\Delta Q_i$ data from a single random time step in the training domain. We thus sample from a $d$-dimensional distribution, preserving correlation between the different $\Delta Q_i$ values at a given time step.
    \item {\bf Independent resampling}: here we resample all marginals independently by selecting a new random time step for each of the $\Delta Q_i$.
    \item {\bf Multivariate Gaussian sampling}: we sample from a multivariate Gaussian distribution (MVG) which was fitted to the joint distribution of the corrections $\Delta Q_i$.
\end{enumerate}

These three sampling approaches are tested for different lengths of the training domain.
All surrogate models are applied to predict the 5000-day QoI distributions, starting from the same initial condition as the reference. Since the surrogate models are stochastic, we perform 10 ``replica simulations", simply executing each LF simulation 10 times without modification.

\begin{figure}[h!]
    \centering
    \includegraphics[width = 0.99 \textwidth]{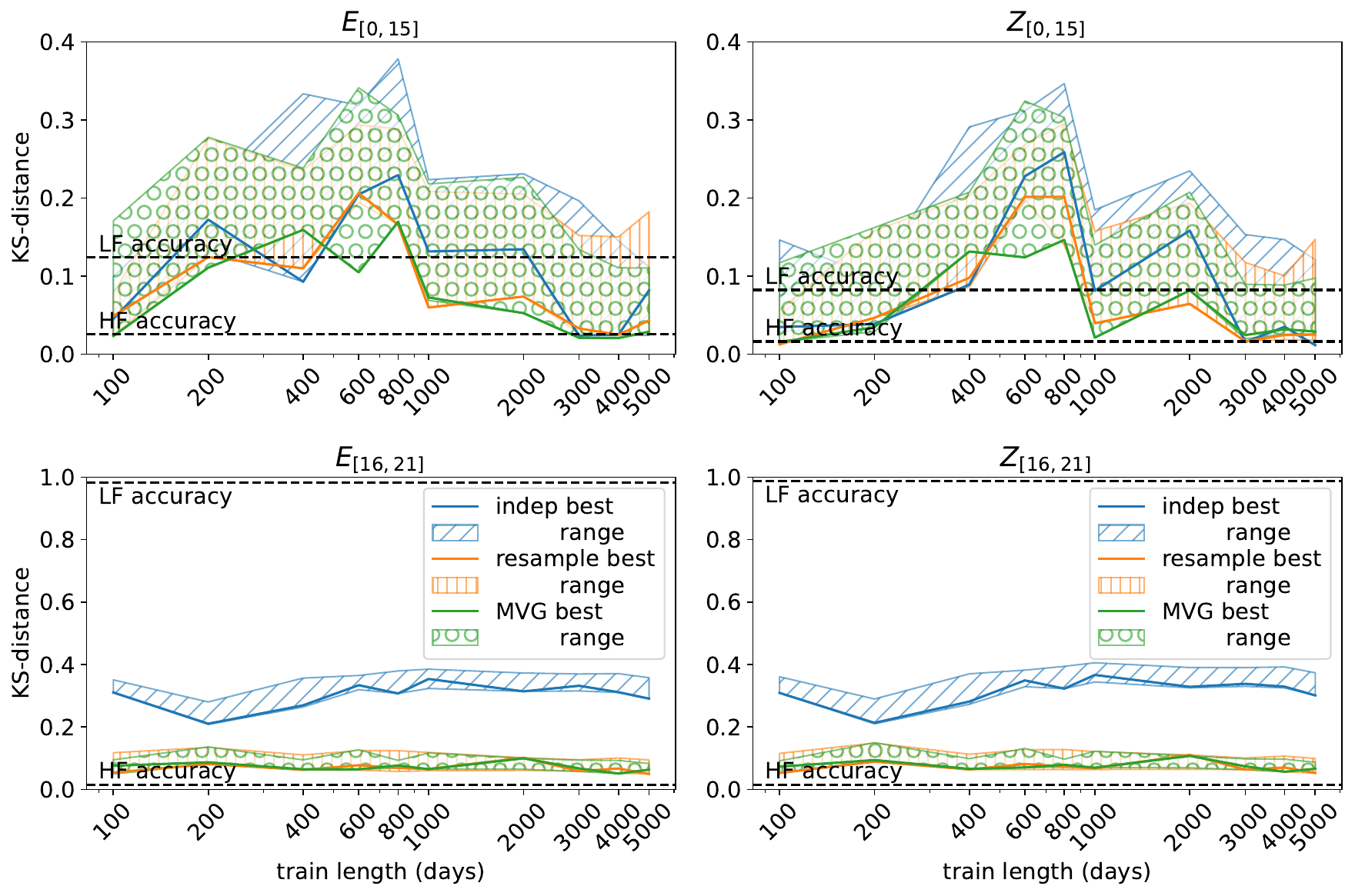}
    \caption{Predictive quality of QoI distributions for the 3 different TO surrogates as a function of the training data size in days. Each surrogate is tested using 10 replicas over a 5000-day simulation period. For each surrogate the ``best" run is plotted, and the shaded area indicates the envelope of the KS-distances over the replicas.}
    \label{fig: KS-distances TO stochastic online}
\end{figure}

The results for the 3 surrogates are visualized in Figure \ref{fig: KS-distances TO stochastic online}, which shows the envelope of the replica KS-distances versus the training data size, alongside the best run of each parameterization. Here ``best" refers to the run which gives the lowest sum over the KS-distances of all QoIs.
``LF accuracy" refers to an LF simulation with no SGS model. The ``HF accuracy" is the KS-distance between the 5000-day QoI distributions in the HF simulation and the full 20000-day reference QoI distributions.

On the large-scale QoIs ($E_{[0,15]}$ and $Z_{[0, 15]}$) the quality of the parameterizations based on at least 3000 training days varies between the LF accuracy and HF accuracy, suggesting that the stochastic surrogates improve the LF simulation most of the time. Such an LF simulation with $\tilde{r}=0$ is already fairly accurate for these 2 QoIs, so any improvement to be had by SGS modelling is relatively small. On the other hand, the parameterizations lead to much improvement on the small-scale QoIs $E_{[16,21]}$ and $Z_{[16, 21]}$, where the KS-distance between the HF simulation and the LF simulation with $\tilde{r}=0$ was almost 1. Here the independent resampling strategy is clearly outperformed by the other approaches.

We see a remarkable pattern between 100 and 800 training days, where the performance deteriorates for our large-scale QoIs ($\Delta E_{[0,15]}$ and $\Delta Z_{[0,15]}$), when including more data. Referring to Figure \ref{fig: convergence of mean dQ}, we can see that the convergence of the relative error in the mean of these corrections lags behind to the corresponding convergence in the small-scale mean of $\Delta E_{[16,21]}$ and $\Delta Z_{[16,21]}$. This indicates that, even though visual inspection of the PDFs does not reveal large discrepancies (see e.g.\ Figure \ref{fig: convergence of dQ and Q_HF}), the TO method can be sensitive to small errors in the mean of the $\Delta Q_i$ noise model. For ease of interpretation the results of Figure \ref{fig: KS-distances TO stochastic online} are also displayed using kernel density estimates in Figure \ref{fig: long term pdfs TO and CNN}, for the simulations with MVG sampling based on 3000 days of training data.

\subsection{CNN results}

\subsubsection{Training error}

The CNNs are trained on snapshots of the coarse-grained solution fields ($\bar{\omega}$, $\bar{\psi}$ and $\bar{r}$) of the same reference simulation as in the previous sections. One snapshot per day (100 HF time steps) is used, as this limits the amount of training data whilst keeping it informative. The training data sets are enriched by including all 90-degree rotations of the snapshots \cite{guan_learning_2023}.

The CNN parameterization is deterministic, yet the training of the CNNs is non-deterministic due to random initialization of the network and random mini-batch selection in the training data. We therefore train 10 replica CNNs to quantify the uncertainty due to stochasticity in training.

\begin{figure}[htp]
    \centering
    \includegraphics[width = 0.6 \textwidth]{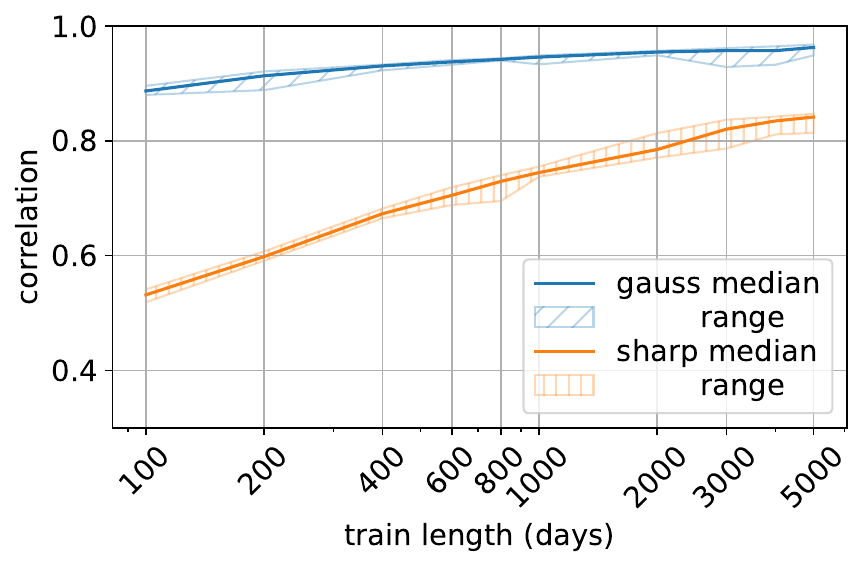}
    \caption{Offline performance of CNNs on withheld test data, as a function of the training data size.}
    \label{fig: offline cnn performance}
\end{figure}
Each trained CNN is evaluated on a withheld test set of snapshots from the HF simulation comprising the data from days 3000 - 3300. Figure \ref{fig: offline cnn performance} shows the test Pearson correlation coefficient between the predicted subgrid scale term and the true subgrid scale term, as a function of the training-data size. As discussed in Section \ref{section: CNN}, we explore both a sharp spectral filter and a Gaussian filter.

The CNNs predict the subgrid scale term from the simulations with the Gaussian filter much better, as expected. The spread in the performance due to the randomness is small. Finally, note that we train for a fixed number of epochs (100, see Section \ref{section: CNN}), so adding more data in the form of extra snapshots does increase the training time.

\subsubsection{CNN prediction results}

\begin{figure}[h!]
    \centering
    \includegraphics[width = 0.9 \textwidth]{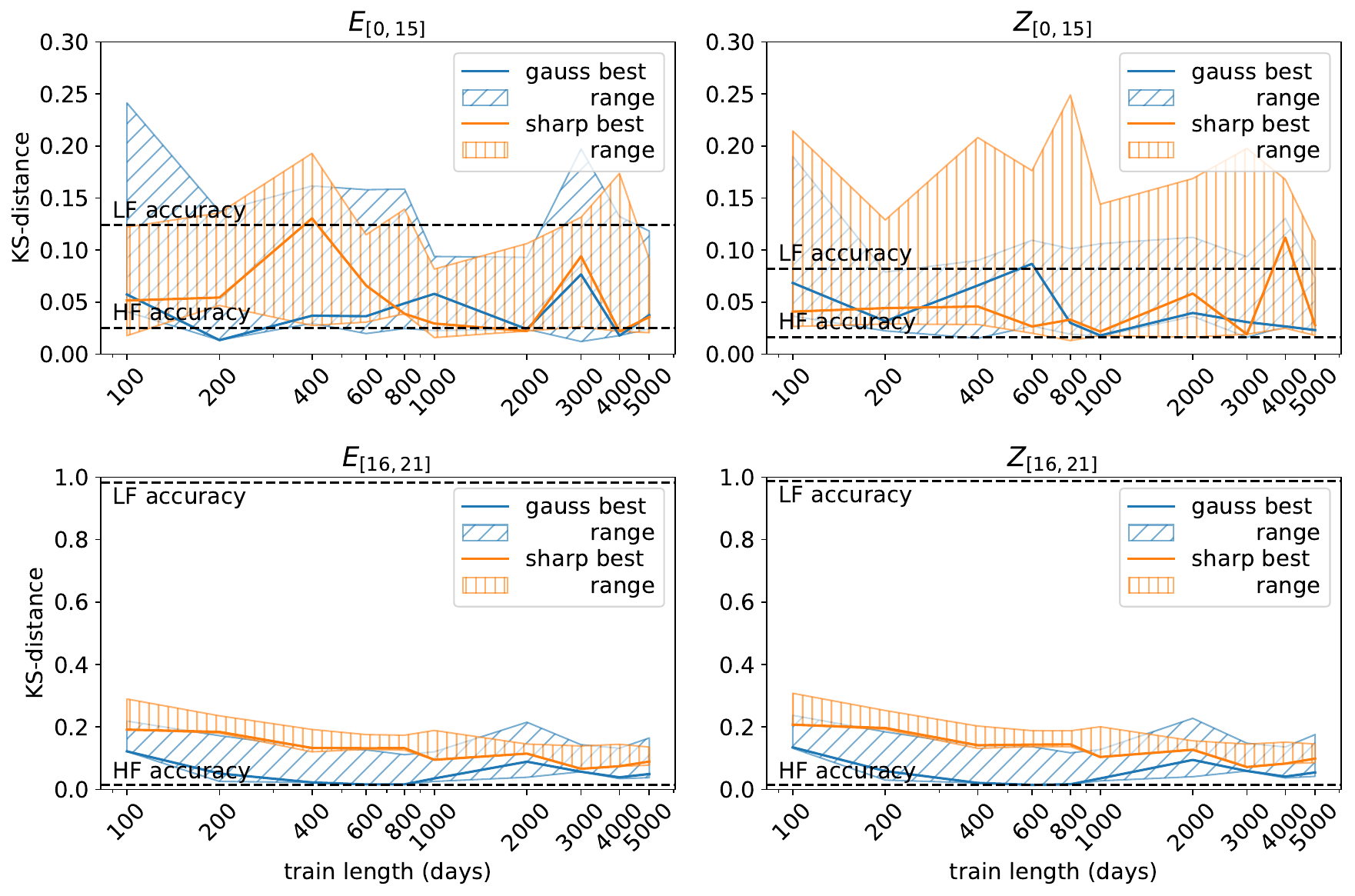}
    \caption{Predictive quality of QoI distributions for the CNN surrogates as a function of the training data size in days. Each surrogate is tested using 10 replicas over a 5000-day simulation period. For each surrogate the ``best" run is plotted, and the shaded area indicates the envelope of the KS-distances over the replicas.}
    \label{fig: KS-distances CNN online}
\end{figure}

The predictive performance of the CNN SGS surrogates will again be measured in terms of how well the long-term distributions of the QoIs are reproduced. Figure \ref{fig: KS-distances CNN online} shows the KS-distance between the 5000-day marginals of the QoIs and the marginals from the 20000-day reference simulation. The spread in the performance due to randomness in training is now much larger compared to Figure \ref{fig: offline cnn performance}. In general, the CNNs with Gaussian filter work better. The only advantage with the sharp filter lies in its smaller range on the small-scale QoIs when trained on large datasets. 

The CNNs with the sharp spectral filter sometimes have notably bad predictive performance for $Z_{[0, 15]}$ despite being trained on large training data sets. To further illustrate the difference between the good and bad performing CNN replicas, Figure \ref{fig: appendix pdfs best and worst CNN 3000 sharp} in \ref{appendix: pdfs best and worst CNN 3000 sharp} contains the kernel density estimates of the best and worst CNN instance trained on 3000 days with the sharp filter. While our simulation does not become unstable as observed elsewhere \cite{beck2019deep, brenowitz2019spatially}, a CNN can cause clear bias in the $Z_{[0, 15]}$ distribution. We emphasize that the large variability in $Z_{[0, 15]}$ performance takes place within replica simulations, such that simply retraining the same CNN can cause significant differences in the predicted outcome. This shows potential unpredictability of {\it a-priori} machine-learning SGS surrogates when coupled to a physics-based partial differential equation.

\subsection{Smagorinsky results}

For a fair comparison, we also optimize the Smagorinsky model in an \textit{a-posteriori} fashion, by tuning $C_s$ to give the best possible long-term distributions for our QoIs over 5000 days in a low-fidelity simulation. For the subgrid scale characteristic length in the $65 \times 65$ low fidelity simulation, we take $\Delta = \frac{2 \pi}{65}$.

\begin{figure}[h!]
    \centering
    \includegraphics[width = 0.9 \textwidth]{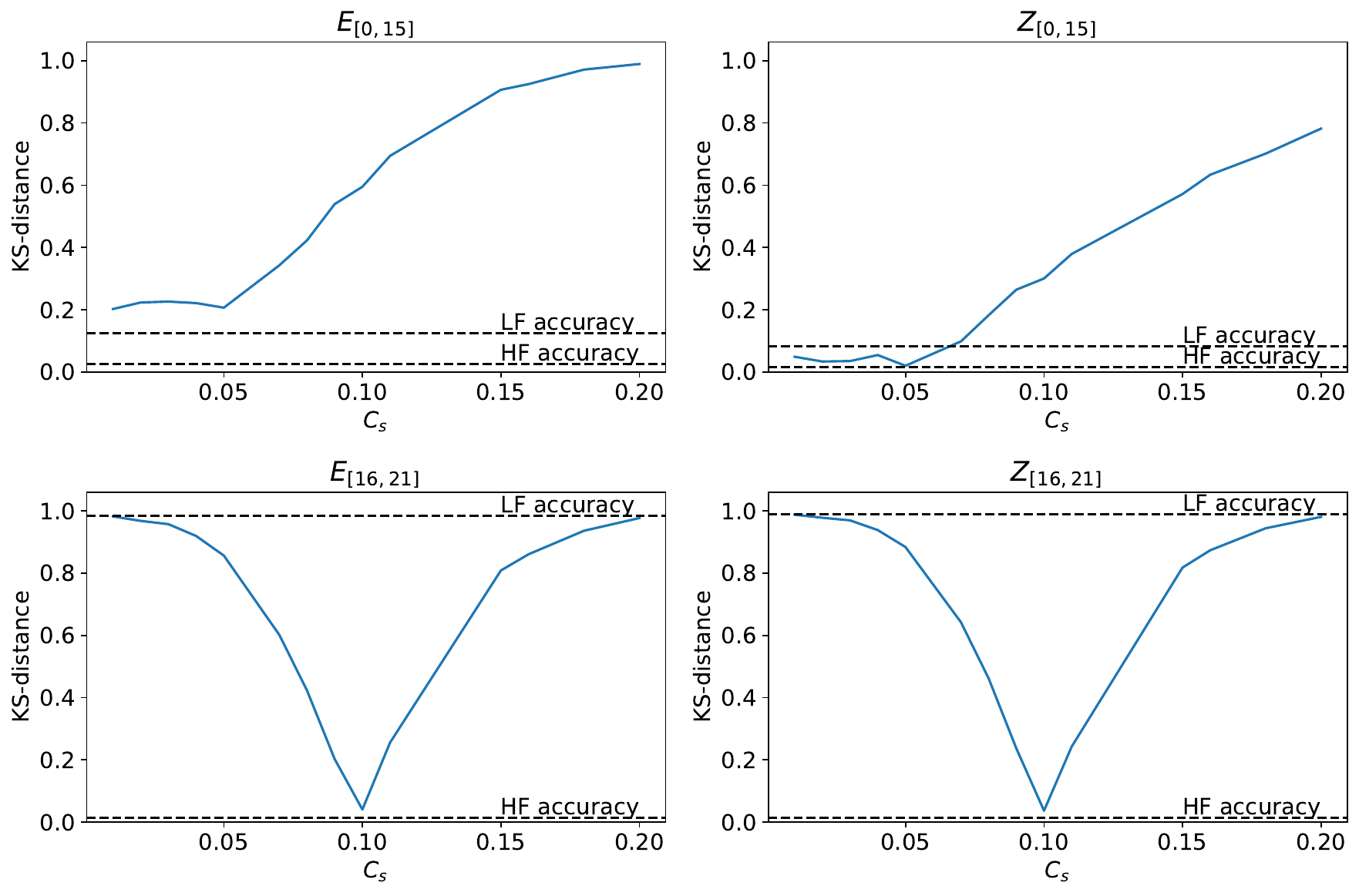}
    \caption{Predictive quality of QoI distributions for the Smagorinsky surrogate as a function of $C_s$.}
    \label{fig: online KS distance smag}
\end{figure}

Figure \ref{fig: online KS distance smag} shows the KS-distances to the reference distributions for the long-term distributions of the QoIs, as a function of $C_s$. As in previous sections these distances are based on 5000-day LF runs and compared to a 20000-day reference simulation. The distributions of the energy and enstrophy in the $[16,21]$ wave number range can almost completely be corrected by using  $C_s = 0.1$. The distribution of the enstrophy in the large scales can also significantly be improved, although this requires $0.02 \leq C_s \leq 0.05$. No $C_s$ value improves upon the LF accuracy with regards to $E_{[0, 15]}$. Given these disjoint optimal parameter values we consider $C_s = 0.1$, which gives the lowest overall summed KS-distance of $0.972$.

\subsection{Comparison of subgrid scale surrogates}

In the previous sections, different versions of three subgrid surrogates were tested. In this section, we will compare the best performing versions of each of these parameterizations in more detail. These are
\begin{enumerate}
    \item The MVG $\Delta Q_i$ surrogate with 3000 days of training data.
    \item The CNN with Gaussian filter trained on 1000 $\bar{r}$ snapshots and their rotations.
    \item The Smagorinsky model with $C_s = 0.1$.
\end{enumerate}

\begin{figure}[h]
\centering
    \includegraphics[width=0.4\textwidth]{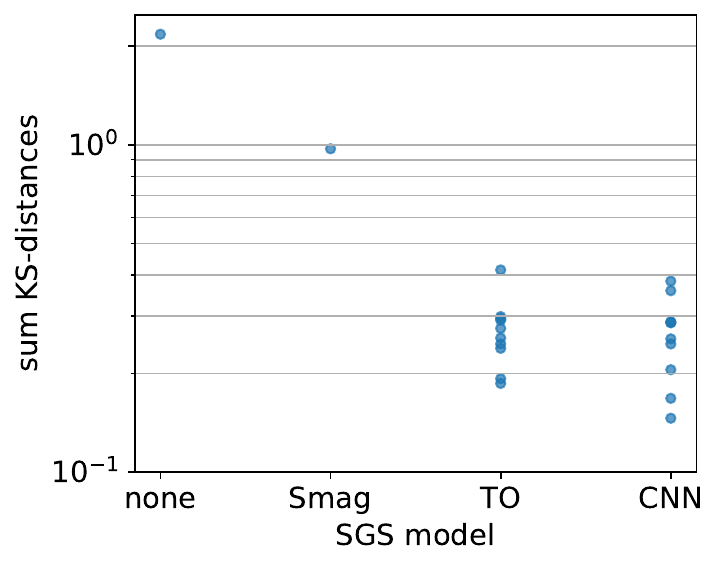}
    \caption{Comparison of summed KS-distances for best versions of SGS surrogates. For surrogates with a stochastic component all replicas are plotted.}
    \label{fig: compare summed KS-dist}
\end{figure}
A comparison of the summed KS-distances for the subgrid surrogates is given in Figure \ref{fig: compare summed KS-dist}. The TO parameterization and CNN parameterization clearly outperform the Smagorinsky parameterization. 

\begin{figure}[htp]
    \centering
    \begin{subfigure}[b]{\textwidth}
        \centering
         \includegraphics[width=\textwidth]{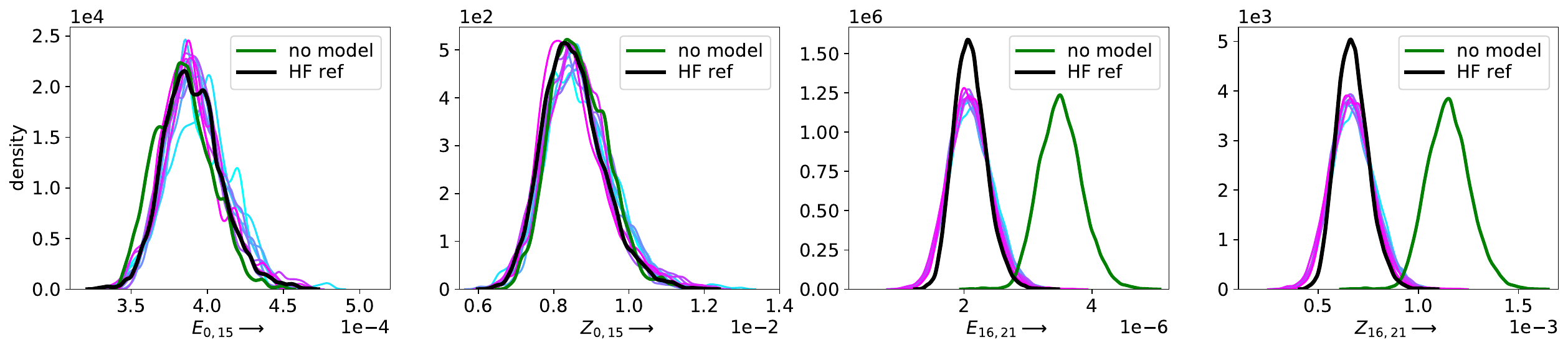}
         \caption{Long-term distributions for TO parameterization with MVG sampling trained on 3000 days.}
         \label{fig: pdfs TO}
    \end{subfigure}
    \begin{subfigure}[b]{1\textwidth}
        \centering
         \includegraphics[width=\textwidth]{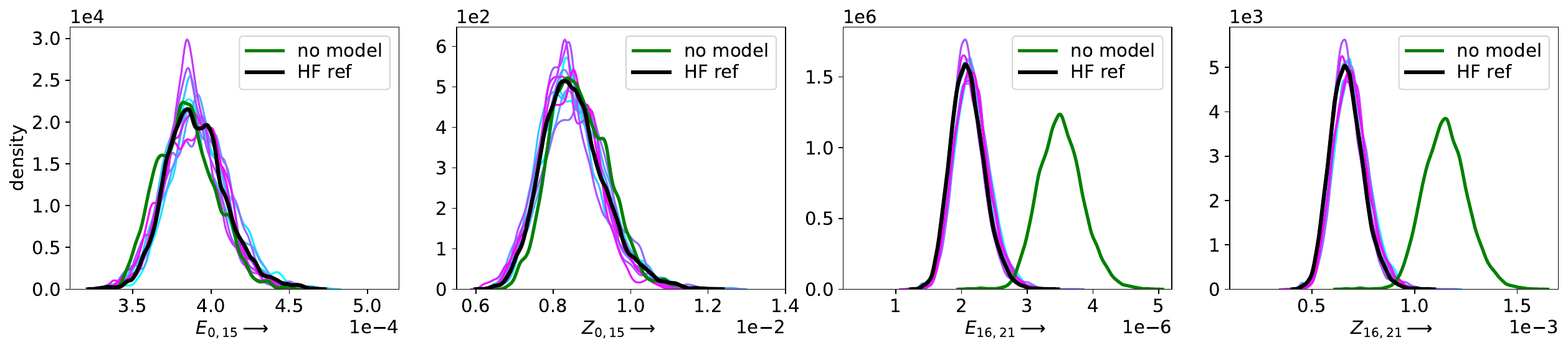}
         \caption{Long-term distributions for CNN parameterization with Gaussian filter trained on 1000 days.}
         \label{fig: pdfs CNN}
    \end{subfigure}
    \caption{Long-term probability-density functions of the 4 QoI.}
    \label{fig: long term pdfs TO and CNN}
\end{figure}

Figure \ref{fig: long term pdfs TO and CNN} shows the 10 long-term replica distributions obtained from the TO and CNN parameterizations, as well as the $\tilde{r}=0$ predictions. Clearly, without an SGS model the statistics for the smaller scales $k\in [16,21]$ are incorrect.
The distributions of the small-scale energy and enstrophy obtained from the TO method are slightly too wide. These distributions have the correct variance when using the CNN parameterization.

\begin{figure}[!tbp]
  \centering
  \begin{minipage}[t]{0.49\textwidth}
    \includegraphics[width=\textwidth]{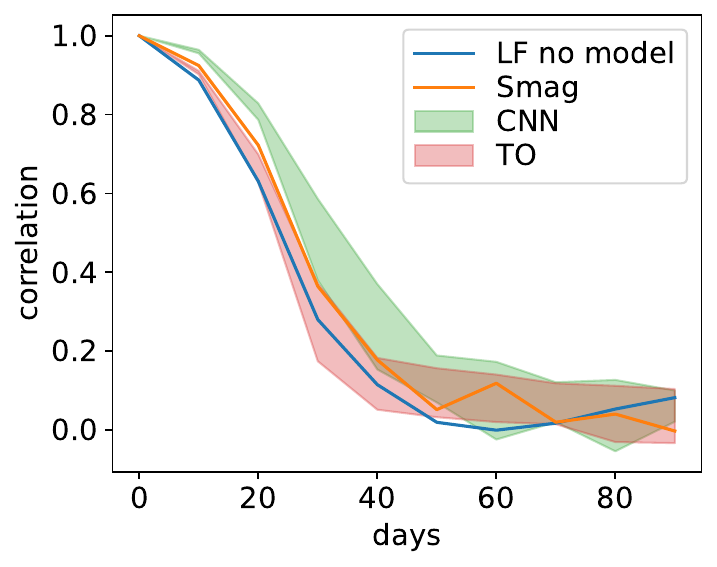}
    \caption{Correlation between LF and HF vorticity.}
    \label{fig: correlation with HF short term}
  \end{minipage}
    \hfill
  \begin{minipage}[t]{0.49\textwidth}
    \includegraphics[width=\textwidth]{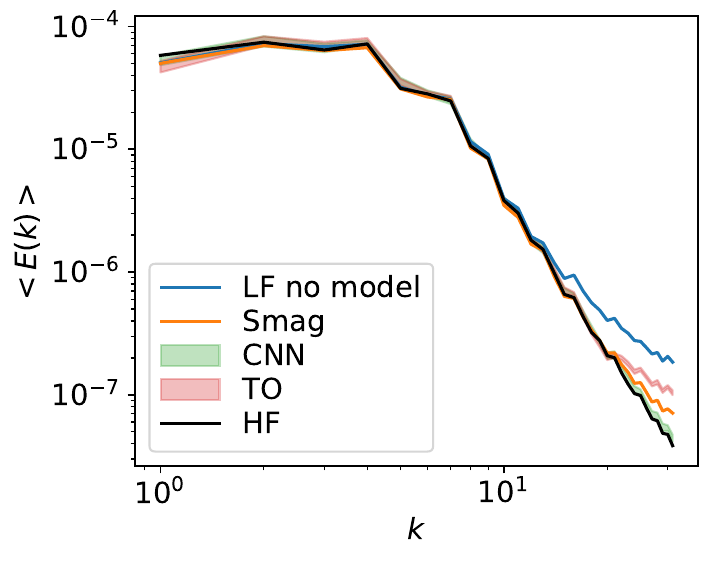}
    \caption{Energy spectra averaged over time.}
    \label{fig: averaged energy spectra}
  \end{minipage}
\end{figure}

\begin{figure}
    \centering
    \includegraphics[width= \textwidth]{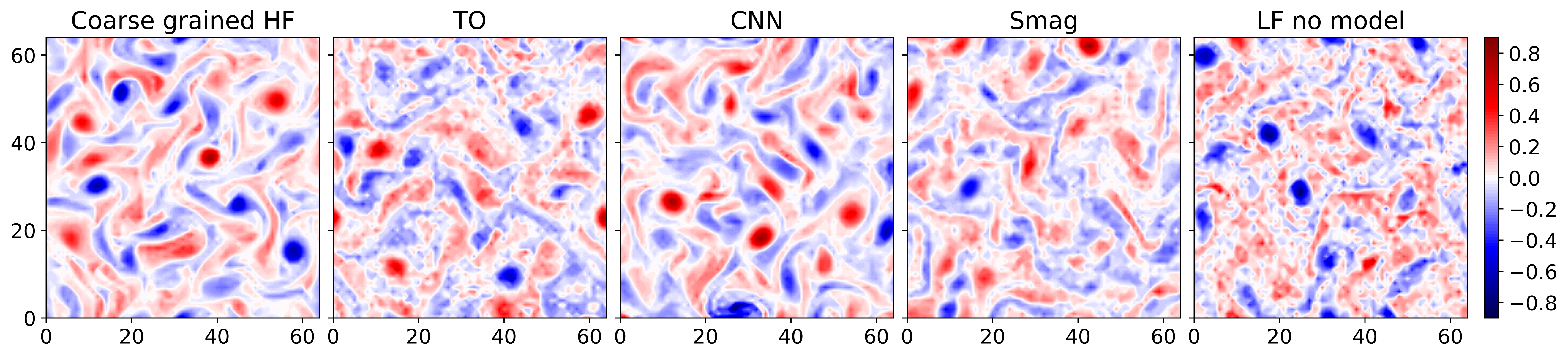}
    \caption{Vorticity fields at day 5000.}
    \label{fig: final vorticity fields}
\end{figure}

It is illustrative to also examine the short-term prediction skill of the subgrid parameterizations. Figure \ref{fig: correlation with HF short term} shows how the LF simulations with different SGS surrogates decorrelate with the HF vorticity solution over time. For each LF simulation every tenth day the correlation coefficient between the solution and the coarse-grained high-fidelity solution is computed. The surrogates with stochastic component are again plotted via their envelope. The decay of correlation with the CNN surrogate is less steep compared to the other surrogates. Still, the TO subgrid parameterization does not decorrelate faster than the LF simulation without SGS term, despite the random error corrections which are added. 

In Figure \ref{fig: averaged energy spectra} the time-averaged energy spectra of the simulations are plotted. These averages are based on 40 snapshots of the solution evenly spaced between days 1000 and 5000 of the simulations. The CNN and TO parameterizations are not averaged over their replica runs. We still plot their envelopes, although time averaging significantly reduced variance between replica simulations. All surrogates (except $\tilde{r}=0$) successfully correct the build-up of energy in the small scales up to approximately $k=20$, for higher wave numbers the CNN parameterizations perform better than the TO parameterizations and the Smagorinsky parameterization. That said, note that the tau-orthogonal surrogates were trained to track QoIs defined on scales up to $\lVert \boldsymbol{k} \rVert = 21$. With this method, one has to decide beforehand what must be learned from data, and guarantees outside the set of predefined QoIs cannot be given.

Finally, one might be interested in visually inspecting the vorticity fields at the end of the simulations. In Figure \ref{fig: final vorticity fields} the vorticity fields at day 5000 are plotted. The $\tilde{r}=0$ simulation leads to a rather non-smooth field. The Smagorinsky and TO parameterizations have smoother fields. The CNN parameterization nicely retains the smoothness of the coarse-grained HF solution. To allow for a fair comparison we applied the inverse of the Gaussian filter to the result of the simulation with the CNN. 

\subsubsection{Computational cost}\label{sec: computational cost}

Coarse-graining is used to alleviate the computational burden of high-accuracy LES. The most costly operations in our pseudospectral solver are (inverse) Fast Fourier transforms \cite{golub2013matrix}, with complexity $O((\frac{3}{2}N)^2 \log((\frac{3}{2}N)^2)$. Theoretically, we can therefore expect the LF solver to compute a time step 20.8 times faster than the HF solver, based on $N_{HF} = 257$, $N_{LF} = 65$. Moreover, the LF solver uses a 10 times larger time step, giving us a total speedup factor of 208. However, adding an SGS surrogate model to the LF simulation might cause serious overhead.

In modern-day computing, the speedup of the solver does not only depend on computational complexity. Multiplications of large arrays can be computed efficiently on a GPU. For smaller problems, communication of data starts taking up relatively more time, which can adversely affect the speedup. In this section, some timing results are presented. Most numerical calculations are done using the PyTorch library \cite{paszke2019pytorch}, which allows for easy offloading of computations to a GPU \footnote{Tests are performed on a NVIDIA GeForce GTX 1070 GPU.}. 

In Table \ref{tab: timing 100 day sim} the runtime of the HF solver and the LF solver with different SGS surrogates are compared for a 100-day simulation. The initialization of the solver is excluded. We also turn off any writing to files.
\begin{table}[ht]
    \centering
    \begin{tabular}{c|cc}
        solver & time [s]\\
        \toprule
        HF & 124\\
        \hdashline
        No model & 1.9 \\
        Smag & 3.4 \\
        CNN & 13.2 \\
        TO resample & 2.8 \\
        TO track & 3.2 \\\bottomrule
    \end{tabular}
    \caption{Time spent on a 100-day simulation.}
    \label{tab: timing 100 day sim}
\end{table}

The LF solver without closure model reaches a 65 times speedup with respect to the HF solver in practice. The Smagorinsky and TO model have an acceptably low overhead. ``TO track" refers to the setup where the TO parameterization is used to track 4 QoIs from a reference simulation. This setup is slightly slower than the online variant with resampling, because it needs to read the reference data from disk.

The CNN subgrid parameterization leads to serious overhead. In fact, one LF time step with CNN is slower than one HF time step, such that all the speedup can be attributed to the coarser time step. Normally, making predictions in batches improves efficiency of neural networks. However, the input of the CNN in the subgrid parametrization depends on the output of the CNN in the previous time step. Due to this recurrent nature of predicting the subgrid term, it is not possible to predict the SGS term for multiple time steps in one batch.

\begin{table}[htp]
    \centering
    \begin{tabular}{cc|c}
        surrogate & train domain [days] & time [s] \\
        \toprule
        CNN & 1000 & 2200\\
        TO track & 3000 & 92.5 \\\bottomrule
    \end{tabular}
    \caption{Training time for data-driven parameterizations.}
    \label{tab: training time}
\end{table}
Before the CNN and TO surrogates can be used, they need to be trained on reference data. Table \ref{tab: training time} reports the training time for a CNN and the TO parameterization. Note that training the TO method involves running an LF simulation that tracks the HF QoI data in the training domain. The TO method needs a longer reference run than the CNN, which leads to extra costs. Taking all this together, the total runtime for finding long term distributions is roughly the same with the TO method and the CNN parametrization. Note that this is problem specific and heavily depends on how much training data both methods need. Also the efficiency of the code implementation and the compute architecture play an important role. Our base solver is an in-house Python code, certainly slower than (more optimized) similar codes written in other languages. Moreover, newer GPUs tend to be optimized for deep learning, leading to faster training and evaluation of the CNNs. The evaluation time of the base solver and the tau-orthogonal method (which also run on the GPU), would likely not profit at the same rate from these optimizations, which would alter the reported relative performance numbers.

To assess the uncertainty in stochastic data-driven surrogates, one might want to run an ensemble of simulations. For the CNN parameterization this would entail performing a new costly training for every ensemble member. For the TO parameterization the online evaluation just needs to be performed with a different random seed.

\section{Conclusion \label{sec:con}}

We presented a study focused on predicting the long-term statistics of global (spatially-integrated) quantities of interest (QoIs) using the tau-orthogonal method in a 2D turbulent simulation. Instead of directly trying to learn the exact unresolved subgrid scale (SGS) term (which is unknown at every point in the computational grid), we replace it by a model with only a handful of unknown degrees of freedom. This small component is the only part which must be learned from (simple and low-dimensional) QoI discrepancy data.

We contrast the tau-orthogonal method with a physics-based Smagorinsky SGS model and a fully data-driven Convolution Neural Network (CNN). To use the tau-orthogonal method outside its training domain, we need a model to extrapolate the time series. We obtained accurate long-term QoI distributions on par with those obtained by the CNN SGS model, by simply modelling the time series as random noise distributed according to the QoI discrepancy distribution, which we estimated in an \textit{a-posteriori} fashion from data. This is despite the fact that the number of unresolved components was compressed by more than a factor of 1000 compared to the CNN. Also, unlike the CNN, the noise model has no tunable parameters. That said, for the simple random-noise models considered, we did require a longer run of the reference model compared to the CNN, in order to capture our large-scale QoI. Our small-scale QoI were already accurately predicted with little data. It would be interesting to examine if using multiple short training trajectories started from cleverly chosen initial conditions can speed up this process. Both the CNN and tau-orthogonal method outperformed the Smagorinsky model in terms of long-term QoI statistics prediction. Moreover, performing low fidelity simulations with the tau-orthogonal method is computationally much cheaper than using the CNN parameterization.

Furthermore, we observed that the CNN could introduce a bias in the long-term statistics, if a sharp spectral filter is used in the resolved-scale equation. Remarkably, this bias could materialize or disappear simply by retraining the CNN surrogate. In a CNN with a Gaussian filter this behavior was less prevalent, possibly because this filter is local in physical space, similar to the convolution filters learned by the CNN. Another possible explanation for this behavior is that the Gaussian filter damps the resolved high wave number contributions, such that the corresponding learning problem is less prone to spectral bias, which is a known issue for CNNs.

A downside of the tau-orthogonal method is that (by construction) it cannot give guarantees about the long-term statistics of quantities outside the predefined set of QoIs. One has to choose beforehand what must be learned from the data. We also observed that with the CNN SGS model, the vorticity fields retain correlation with the high-fidelity reference fields for a longer period. This is mainly because the local spatial correlation of the vorticity fields is better captured when using the CNN SGS model. We expect that this is due to the choice of spatio-temporal basis functions in the tau-orthogonal method, which is the subject of future research.

Thus far, the tau-orthogonal method has only been applied to the 2D barotropic Navier-Stokes equations. While this has relevance for geophysical flow simulations, an essential research direction is the application to 3D flow problems. Such a realistic flow case is important to truly evaluate the promise of the data-driven method. Here, the computational costs and memory requirements for the base solver increase tremendously, due to the cubic scaling of the problem size. This means that the relative decrease in computational cost from coarse graining also is significantly larger than in the two-dimensional case, i.e. scaling cubically instead of quadratically. CNN parameterizations for such problems rely on three-dimensional kernels, drastically increasing the network size and training cost compared to their two-dimensional counter parts. On the other hand, in the tau-orthogonal method the dimension of the learning problem is independent of the spatial dimension of the differential equation under consideration. Performance wise, the only significant extra costs are incurred by taking inner products of three-dimensional vector fields. Whether the tau-orthogonal method is able to accurately capture quantities of interest in a large 3D simulation is the subject of ongoing research.

\section*{Software \& data availability}

The implementation of the methods which are presented in this paper is publicly available at \cite{hoekstra_2024}. This code repository includes a notebook, which takes the user step by step through the experiments in the paper. 

\section*{Acknowledgements}

This research is funded by the Netherlands Organization for Scientific Research (NWO) through the ENW-M1 project ``Learning small closure models for large multiscale problems” (OCENW.M.21.053). This work used the Dutch national e-infrastructure with the support of the SURF Cooperative using grant no.\ EINF-9277. We further made us of the ARCHER2 UK National Supercomputing Service under the UK EPSRC SEAVEA grant (EP/W007762/1).

\section*{CRediT author statement}

{\bf Rik Hoekstra}: Conceptualization, Data curation, Formal analysis, Investigation, Methodology, Resources, Software, Validation, Visualization, original draft, review \& editing. {\bf Daan Crommelin}: Conceptualization; Methodology; Project administration; Supervision, review \& editing. {\bf Wouter Edeling}: Conceptualization; Formal analysis; Funding acquisition; Investigation; Methodology; Project administration;   Supervision; Validation; review \& editing.

\section*{Declaration of competing interest}

The authors declare that they have no known competing financial interests or personal relationships that could have appeared to influence the work reported in this paper.

\section*{Declaration of Generative AI and AI-assisted technologies in the writing process}

The authors declare that no generative AI and/or AI-assisted technologies were used in the writing process.

\bibliographystyle{plainnat}
\bibliography{refs}

\appendix
\section{Problem setup} \label{app: problem setup}
We solve \eqref{eq:HF} and \eqref{eq:HF2} on a double periodic domain of size $2\pi L \times 2\pi L$. Here, $L=6.371\times10^6 [m]$ is chosen as the Earth's radius. Using the inverse of the Earth's angular velocity $\Omega^{-1}$ as a time scale, $2\pi$ non-dimensionalized units of time can be related to a single ``day" as $2\pi\approx24\times60^2\times\Omega$, where $\Omega=7.292\times 10^{-5}[s^{-1}]$. The diffusion and forcing constants, $\mu$ and $\nu$, are chosen such that Fourier modes at the smallest resolved scale are damped with an e-folding time scale of 5 and 90 days respectively: $\nu \approx 4.4 \cdot 10^{-6}$,  $\mu \approx 1.8 \cdot 10^{-3}$, see \cite{edeling_reducing_2020,verkley_maximum_2016} for more detail. The forcing is applied to a single Fourier mode $F = 2^{\frac{3}{2}}\cos(5x)\cos(5y)$.

For our experiments, initial conditions in the statistically stationary state of the system are used, which are obtained after a 300-day spin-up period. In turn, the spin-up period is initialized at $t=0$ from
\begin{equation}
    \begin{split}
        \omega(x,y;t=0) = & \sin(4x) \sin(4y) + 0.4 \cos(3x)\cos(3y) \\ 
        & + 0.3 \cos(5x) \cos(5y) + 0.02 \sin(x) + 0.02 \cos(y).
    \label{eq:IC}
    \end{split}
\end{equation}

\section{Pseudospectral treatment advection term} \label{app:pseudospectral treatment}
To remove aliasing errors in the non-linear term $J$, see \eqref{eq: J}, we apply the well-known ``3/2 rule”\cite{peyret_vorticity-streamfunction_2002}, which entails the following steps. First, each $N \times N$ field in the Fourier domain is extended with zeros to an $\frac{3}{2}N \times \frac{3}{2} N$ field. The individual derivatives in $J$ are computed in the Fourier domain, after which an inverse Fourier transform $\mathcal{F}^{-1}$ is applied to multiply the various derivatives of $J$ in the physical domain. After that, the result is transformed back to the Fourier domain. Finally, only the smallest $N \times N$ wave numbers are retained, discarding all wave numbers for which we earlier introduced zeros. This procedure adds numerical dissipation since we allow energy to be transferred to the extended wave numbers and then set them to zero.

\section{Analytic tau-orthogonal surrogate form \label{app:analytic}}

Here we derive the analytic expression for $\tilde{r}$. In our case $R_1=R_2=R_{[0,15]}$ are disjoint from $R_3=R_4=R_{[16, 21]}$, see Figures \ref{fig:filters_b}-\ref{fig:filters_c}. Although we solve linear systems like \eqref{eq:c_1j} numerically in practice, disjoint filters do simplify the analytic expression that underpins $\tilde{r}$, as $c_{ij}=0$ for any disjoint $R_i$ and $R_j$. To see this, note that for our particular case, \eqref{eq:c_1j} simplifies to
\small
\begin{align}
    \left[
    \begin{matrix}
        \left(V_2, R_2 T_2\right) & 0 & 0 \\ 
        0 &  \left(V_3, R_3 T_3\right) & \left(V_3, R_3 T_4\right) \\
        0 &  \left(V_4, R_4 T_3\right) & \left(V_4, R_4 T_4\right)
    \end{matrix}\right]
    \left[
    \begin{matrix}
    c_{12} \\ c_{13} \\ c_{14}
    \end{matrix}
    \right] = 
    \left[
    \begin{matrix}
    - \left(V_2, R_2 T_1\right)
    \\ 0 \\ 0
    \end{matrix}
    \right],
     \label{eq:c_1j_2}
\end{align} 
\normalsize
such that $c_{12}$ is decoupled from the other coefficients, and $[c_{13}, c_{14}]^T=[0, 0]^T$ is the trivial null space solution of the $2\times2$ (symmetric) submatrix displayed in \eqref{eq:c_1j_2}, provided it is full rank. In addition to $c_{ii}=1$, the other 4 non-zero solutions $c_{12}, c_{21}, c_{34}$ and $c_{43}$ are given by
\begin{align}
    c_{ij} = -\frac{\left(V_j, R_jV_i\right)}{\left(V_j, R_jV_j\right)},
\end{align}
\noindent
which has the typical form of an orthogonal projection coefficient. With
\begin{align}
    &T_1=V_1 = -\bar{\psi}_{[0, 15]}, \quad
    T_2=V_2 = \bar{\omega}_{[0, 15]}, \nonumber\\
    &T_3=V_3 = -\bar{\psi}_{[16, 21]}, \quad
    T_4=V_4 = \bar{\omega}_{[16, 21]}, \label{eq: T and V}   
\end{align}
\noindent
the final expression for our reduced SGS term becomes
\begin{align}
    \tilde{r} = &-\frac{\D Q^u_1}{\D t} \left(\frac{ Z_{[0, 15]}}{{2}S_{[0, 15]}Z_{[0, 15]} - {2} E^2_{[0, 15]}}\right)\left(\bar\psi_{[0, 15]} + \frac{E_{[0, 15]}}{Z_{[0, 15]}}\bar\omega_{[0, 15]}\right)  \nonumber\\
    &+\frac{\D Q^u_2}{\D t}\left(\frac{S_{[0, 15]}}{{2}S_{[0, 15]}Z_{[0, 15]} - {2}E^2_{[0, 15]}}\right)\left(\bar\omega_{[0, 15]} + \frac{E_{[0, 15]}}{S_{[0, 15]}}\bar\psi_{[0, 15]}\right)  \nonumber\\
    &- \frac{\D Q^u_3}{\D t} \left(\frac{ Z_{[16, 21]}}{{2} S_{[16, 21]}Z_{[16, 21]} - {2} E^2_{[16, 21]}}\right)\left(\bar\psi_{[16, 21]} + \frac{E_{[16, 21]}}{Z_{[16, 21]}}\bar\omega_{[16, 21]}\right) \nonumber\\
    &+ \frac{\D Q^u_4}{\D t} \left(\frac{S_{[16, 21]}}{{2}S_{[16, 21]}Z_{[16, 21]} - 
    {2} E^2_{[16, 21]}}\right)\left(\bar\omega_{[16, 21]} + \frac{E_{[16, 21]}}{S_{[16, 21]}}\bar\psi_{[16, 21]}\right). 
\end{align}
\noindent
Here, we define the squared stream function as $S_{[l, m]}:=\left(\bar\psi_{[l, m]}, \bar\psi_{[l, m]}\right)/2$. 

\section{Error in QoI after TO correction} \label{app:error in correction}

We derive a bound for the difference between the HF reference trajectories and the LF trajectories of the QoIs at time $t_n$ when we track the reference with our predictor-corrector scheme.
Recall that $\boldsymbol{\Delta Q}^{n+1}$ is defined as vector containing the difference between the predicted values of the QoIs and their reference values
\small
\begin{equation}
    \Delta Q_i^{n+1} := Q_i(\bar{\psi}^{n+1}_{HF}, \bar{\omega}^{n+1}_{HF}) - Q_i(\bar{\psi}^{n+1^*}, \bar{\omega}^{n+1^*}).
\end{equation}
\normalsize
Let us first consider the LF trajectory of the enstrophy in the large scales at time $t_{n+1}$:
\small
\begin{align}
    Z_{LF[0,15]}^{n+1} =& \frac{1}{2}\left(R_2 \bar{\omega}^{n+1}_{LF}, R_2 \bar{\omega}^{n+1}_{LF}\right).
\end{align}
\normalsize
Writing this in terms of the predicted fields and the correction term we get
\small
\begin{align}
    Z_{LF[0,15]}^{n+1} =& \frac{1}{2}\left(R_2 \bar{\omega}^{n+1^*}+R_2 \sum_{i=1}^{d} \tau_i P_i^*, \quad R_2 \bar{\omega}^{n+1^*}+R_2 \sum_{i=1}^{d} \tau_i P_i^*\right) \nonumber\\
                     =& \frac{1}{2}\left(R_2 \bar{\omega}^{n+1^*}, R_2 \bar{\omega}^{n+1^*} \right) +
                        \left(R_2 \bar{\omega}^{n+1^*}, R_2 \sum_{i=1}^{d} \tau_i P_i^* \right)
                         +\frac{1}{2}\left(R_2 \sum_{i=1}^{d} \tau_i P_i^*, R_2 \sum_{i=1}^{d} \tau_i P_i^*\right).
\end{align}
\normalsize
Here $P_i^* = P_i(\bar{\omega}^{n+1^*})$.
We substitute \eqref{eq: tau in terms of dq} and introduce $V_i^* = V_i(\bar{\omega}^{n+1^*})$, such that the sum in the second term disappears due to our orthogonality constraints \eqref{eq: orthogonality constr}:
\small
\begin{align}
   Z_{LF[0,15]}^{n+1} =& Q_2(\bar{\psi}^{n+1^*}, \bar{\omega}^{n+1^*}) + 
                        \left(R_2 \bar{\omega}^{n+1^*}, R_2 \frac{\Delta Q_2^{n+1}}{\left( V_2^*, R_2 P_2^* \right)} P_2^* \right) \nonumber\\
                        &+\frac{1}{2}\left(R_2 \sum_{i=1}^{d} \tau_i P_i^*, R_2 \sum_{i=1}^{d} \tau_i P_i^*\right).
\end{align}
\normalsize
Using $V_2^* = R_2\bar{\omega}^{n+1^*}$:
\small
\begin{align}
    Z_{LF[0,15]}^{n+1} =& Q_2(\bar{\psi}^{n+1^*}, \bar{\omega}^{n+1^*}) + 
                        \left(R_2 \bar{\omega}^{n+1^*}, R_2 \frac{\Delta Q_2^{n+1}}{\left( R_2 \bar{\omega}^{n+1^*}, R_2 P_2^* \right)} P_2^* \right) \nonumber \\
                        &+\frac{1}{2}\left(R_2 \sum_{i=1}^{d} \tau_i P_i^*, R_2 \sum_{i=1}^{d} \tau_i P_i^*\right).
\end{align}
\normalsize

Noting that $\Delta Q_2$ and its denominator are scalar quantities, we see that the second term can be simplified drastically
\small
\begin{align}
   Z_{LF[0,15]}^{n+1} =& Q_2(\bar{\psi}^{n+1^*}, \bar{\omega}^{n+1^*}) + 
                        \Delta Q_2^{n+1}
                        +\frac{1}{2}\left(R_2 \sum_{i=1}^{d} \tau_i P_i^*, R_2 \sum_{i=1}^{d} \tau_i P_i^*\right) \nonumber\\
                    =& Q_2(\bar{\psi}^{n+1}_{HF}, \bar{\omega}^{n+1}_{HF}) 
                        +\frac{1}{2}\left(R_2 \sum_{i=1}^{d} \tau_i P_i^*, R_2 \sum_{i=1}^{d} \tau_i P_i^*\right) \nonumber\\
                    =& Z_{HF[0,15]}^{n+1}
                        +\frac{1}{2}\left( \sum_{i=1}^{d} \frac{\Delta Q_i^{n+1}}{\left( V_i^*, R_i P_i^* \right)} R_2 P_i^*,  \sum_{i=1}^{d} \frac{\Delta Q_i^{n+1}}{\left( V_i^*, R_i P_i^* \right)} R_2 P_i^*\right).
\end{align}
\normalsize
Substituting \eqref{eq: T and V}, we rewrite the error term
\small
\begin{align}
    &Z_{LF[0,15]}^{n+1} = (1+\epsilon) Z_{HF[0,15]}^{n+1},\\
&\epsilon= \frac{Z_{[0,15]}^{n+1^*} (\Delta Q_1^{n+1})^2 + S_{[0,15]}^{n+1^*} (\Delta Q_2^{n+1})^2 - 2 E_{[0,15]}^{n+1^*} \Delta Q_1^{n+1}\Delta Q_2^{n+1}}{(R_1\bar{\psi}^{n+1^*},R_1\bar{\psi}^{n+1^*})(R_1\bar{\omega}^{n+1^*},R_1\bar{\omega}^{n+1^*})-(R_1\bar{\psi}^{n+1^*},R_1\bar{\omega}^{n+1^*})^2} \frac{1}{Z^{n+1^*}_{[0,15]}+\Delta Q_2^{n+1}}.\label{eq: error term tussen 1}
\end{align}
\normalsize
Here we use the energy, enstrophy and squared stream function of the predicted fields, denoted by the superscript $n+1^*$.\\
The generalised angle between $\bar{\omega}$ and $\bar{\psi}$ can be defined based on inner products as
\begin{equation}
    \cos{\theta_{\bar{\omega}, \bar{\psi}}} = \frac{(\bar{\psi}, 
\bar{\omega})}{\lVert \bar{\psi} \rVert \lVert \bar{\omega} \rVert}.
\end{equation}
If $-c \leq \cos{\theta_{\bar{\omega}, \bar{\psi}}} \leq c$ for $c \in [0,1)$, the first part of the denominator in \eqref{eq: error term tussen 1} does not equal 0 and the fields are not linearly dependent.

We get the following bound on the error term
\begin{equation}
\lvert \epsilon \rvert \leq \bigg\lvert \frac{Z_{[0,15]}^{n+1^*} (\Delta Q_1^{n+1})^2 + S_{[0,15]}^{n+1^*} (\Delta Q_2^{n+1})^2 - 2 E_{[0,15]}^{n+1^*} \Delta Q_1^{n+1}\Delta Q_2^{n+1}}{4(1-c^2) S_{[0,15]}^{n+1^*} Z_{[0,15]}^{n+1^*}(Z^{n+1^*}_{[0,15]}+\Delta Q_2^{n+1})} \bigg\rvert.
\end{equation}
We drop the $\Delta Q_2^{n+1}$ term from the denominator and use the Cauchy-Schwarz inequality to bound $E$ in terms of $S$ and $Z$:
\begin{equation}
\lvert \epsilon \rvert \leq 
\frac{1}{4(1-c^2)} \bigg( \bigg\lvert \frac{(\Delta Q_1^{n+1})^2}{Z_{[0,15]}^{n+1^*} S_{[0,15]}^{n+1^*}} \bigg\rvert  +  \bigg\lvert \frac{(\Delta Q_2^{n+1})^2}{Z_{[0,15]}^{n+1^*} Z_{[0,15]}^{n+1^*}} \bigg\rvert +  2\bigg\lvert \frac{\Delta Q_1^{n+1}\Delta Q_2^{n+1}}{Z_{[0,15]}^{n+1^*} E_{[0,15]}^{n+1^*}}  \bigg\rvert \bigg).
\end{equation}
This bound gives clear requirements for accurate tracking of our chosen QoIs. Firstly, the fields $R_1\bar{\psi}^{n+1^*}$ and $R_1\bar{\omega}^{n+1^*}$ should not be close to linearly dependent. Moreover, the size of the corrections $\Delta Q_1$ and $\Delta Q_2$ should be small compared to the predicted values for $Z_{[0,15]}$, $S_{[0,15]}$ and $E_{[0,15]}$.

For the LF trajectory of the energy in the large scales $E_{[0,15]}$ we have
\small
\begin{align}
    E_{LF[0,15]}^{n+1} =& -\frac{1}{2}\left(R_1 \bar{\omega}^{n+1}, R_1 \bar{\psi}^{n+1}\right) = -\frac{1}{2}\left(R_1 \nabla^2 \bar{\psi}^{n+1}, R_1 \bar{\psi}^{n+1}\right) \nonumber\\
                     =& -\frac{1}{2}\left(R_1 \nabla^2 (\bar{\psi}^{n+1^*}+C), R_1 (\bar{\psi}^{n+1^*}+C) \right),
\end{align}
\normalsize
where $C = \nabla^{-2} \sum_{i=1}^{d} \tau_i P_i^*$.

We define our discretized Laplacian and its inverse on a function $f:\mathbb{R}^2 \to \mathbb{R}$ and its discrete Fourier expansion $\hat{f}(\boldsymbol{k})$ as $\nabla^2 f = \mathcal{F}^{-1}(\hat{\nabla}^2 \hat{f}(\boldsymbol{k}))$ and $\nabla^{-2} f = \mathcal{F}^{-1}(\hat{\nabla}^{-2} \hat{f}(\boldsymbol{k}))$, where
\small
\begin{align}
    \hat{\nabla}^2 \hat{f}_{\boldsymbol{k}} &= -\lVert \boldsymbol{k} \rVert_2^2 \hat{f}_{\boldsymbol{k}},\\
    \hat{\nabla}^{-2} \hat{f}_{\boldsymbol{k}}&= \begin{cases}
        0, &\text{ if } \boldsymbol{k} = \boldsymbol{0},\\
        -\frac{1}{\lVert \boldsymbol{k} \rVert_2^2} \hat{f}_{\boldsymbol{k}}, &\text{ else},
    \end{cases}
\end{align}
\normalsize
in line with our discretization of the governing equations \cite{peyret_vorticity-streamfunction_2002}. Note that these discrete operators retain the self-adjoined property of the Laplacian under the approximation of the inner-products in \eqref{eq:EZ}. We expand the inner-product, use the self-adjoinedness, followed by the same steps we took in the enstrophy derivation:
\small
\begin{align}
    E_{LF[0,15]}^{n+1} =& -\frac{1}{2}\left(R_1 \nabla^2 \bar{\psi}^{n+1^*}, R_1 \bar{\psi}^{n+1^*} \right) -
                        \left(R_1 \bar{\psi}^{n+1^*}, R_1 \nabla^2 C \right)
                         -\frac{1}{2}\left( R_1 \nabla^2 C, R_1 C\right) \nonumber\\
                     =& Q_1(\bar{\psi}^{n+1^*}, \bar{\omega}^{n+1^*})
                        -\left(R_1 \bar{\psi}^{n+1^*}, R_1 \sum_{i=1}^{d} \frac{\Delta Q_i^{n+1}}{\left( V_i^*, R_i P_i^* \right)} P_i^* \right) \nonumber\\
                        & -\frac{1}{2}\left( R_1 \nabla^2 C, R_1 C\right) \nonumber\\
                     =& Q_1(\bar{\psi}^{n+1^*}, \bar{\omega}^{n+1^*}) - 
                        \left(R_1 \bar{\psi}^{n+1^*}, R_1 \frac{\Delta Q_1^{n+1}}{\left( -R_1 \bar{\psi}^{n+1^*}, R_1 P_1^* \right)} P_1^* \right) \nonumber\\
                        & -\frac{1}{2}\left( R_1 \nabla^2 C, R_1 C\right) \nonumber\\
                    =& E_{HF[0,15]}^{n+1}
                         -\frac{1}{2}\left( R_1 \nabla^2 C, R_1 C\right).
\end{align}
\normalsize
Finally we substitute $C$ to obtain
\small
\begin{align}
   E_{LF[0,15]}^{n+1} =& E_{HF[0,15]}^{n+1}
                         -\frac{1}{2}\left(  \sum_{i=1}^{d} \tau_i R_1 P_i^*, \nabla^{-2} \sum_{i=1}^{d} \tau_i R_1 P_i^*\right).
\end{align}
\normalsize
We rewrite the error as a relative error
\begin{equation}
    E_{LF[0,15]}^{n+1} = (1+\epsilon) E_{HF[0,15]}^{n+1}.
\end{equation}
Using the same bounds as before we get
\begin{align}
    \lvert \epsilon \rvert \leq & 
    \frac{1}{4(1-c^2)^2} \bigg(
    \bigg\lvert \frac{(\Delta Q_1^{n+1})^2}{Z_{[0,15]}^{n+1^*} S_{[0,15]}^{n+1^*}} \bigg\rvert  + 
    \bigg\lvert \frac{(\Delta Q_2^{n+1})^2}{Z_{[0,15]}^{n+1^*} Z_{[0,15]}^{n+1^*}} \bigg\rvert + 
    2\bigg\lvert \frac{\Delta Q_1^{n+1}\Delta Q_2^{n+1}}{Z_{[0,15]}^{n+1^*} E_{[0,15]}^{n+1^*}}  \bigg\rvert \\ \nonumber
    & + \bigg\lvert \frac{(\Delta Q_1^{n+1})^2 F}{E_{[0,15]}^{n+1^*} (S_{[0,15]}^{n+1^*})^2} \bigg\rvert 
    +  \bigg\lvert \frac{(\Delta Q_2^{n+1})^2 F}{Z_{[0,15]}^{n+1^*} E_{[0,15]}^{n+1^*} S_{[0,15]}^{n+1^*}} \bigg\rvert + 
    2\bigg\lvert \frac{\Delta Q_1^{n+1}\Delta Q_2^{n+1} F}{Z_{[0,15]}^{n+1^*} (S_{[0,15]}^{n+1^*})^2}  \bigg\rvert \bigg),
\end{align}
where $F = (R_1 \bar{\psi}^{n+1^*}, R_1 \nabla^{-2} \bar{\psi}^{n+1^*})/2$. Note that the first three terms in this bound are the same as in the enstrophy bound. The factor $F$ in the additional terms is well-behaved, due to the strong smoothing effect of the inverse Laplacian. We can bound $F$ by
\begin{equation}
    \lvert F \rvert = \lvert \frac{1}{2} (R_1 \bar{\psi}^{n+1^*}, R_1 \nabla^{-2} \bar{\psi}^{n+1^*}) \rvert \leq \frac{1}{2 k_{min}^2} (R_1 \bar{\psi}^{n+1^*}, R_1 \bar{\psi}^{n+1^*}) = \frac{1}{k_{min}^2} S_{[0,15]}^{n+1^*},
\end{equation}
where $k_{min}=1$ is the length of the wave number vector corresponding to the largest possible wave in $R_1 \bar{\psi}^{n+1^*}$.

For the energy and enstrophy in the small scales, the accuracy can be derived analogously, resulting in the same bounds with different filters and $k_{min}=16$.

\section{PDFs of long term distributions CNN for sharp filter} \label{appendix: pdfs best and worst CNN 3000 sharp}

\begin{figure}[htp]
    \centering
    \includegraphics[width = \textwidth]{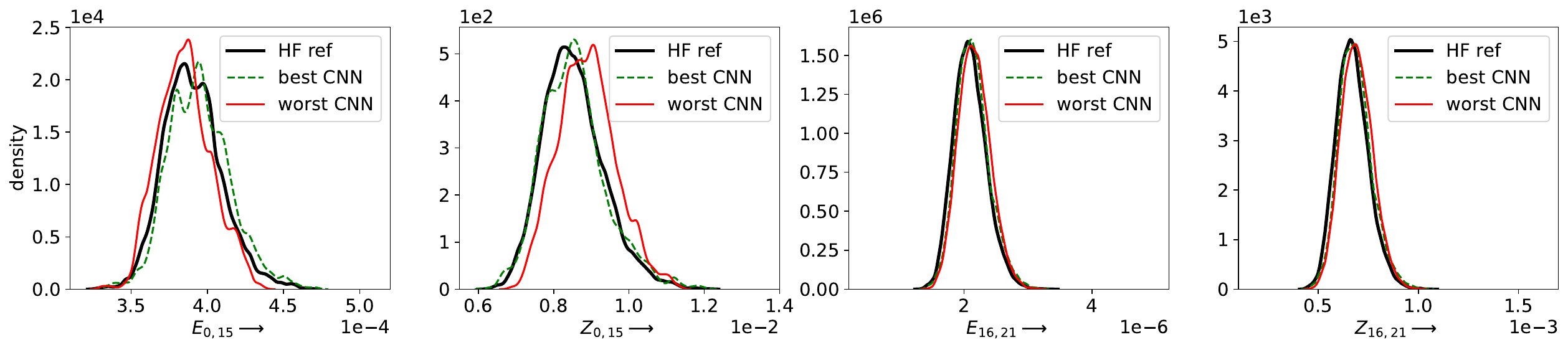}
    \caption{Long-term distributions for both the best and worst CNN parameterization with a sharp spectral filter, trained on 3000 days. Best and worst CNNs are selected based on the summed KS-distance.}
    \label{fig: appendix pdfs best and worst CNN 3000 sharp}
\end{figure}

Here we display the best and worst performing CNN surrogates with a sharp spectral filter trained on 3000 days. Note that a bias in $Z_{[0,15]}$ can develop, see Figure \ref{fig: appendix pdfs best and worst CNN 3000 sharp}.

\end{document}